% plain TeX file of paper ``q-Special functions, a tutorial''
% by Tom H. Koornwinder
%
\magnification 1200
\vsize 9.1truein

\input epsf

\input amssym.def

\overfullrule0pt

\def\ifundefined#1{\expandafter\ifx\csname#1\endcsname\relax}

\newtoks\sectionnumber
\newcount\equationnumber
\newcount\thnumber
\newcount\refnum
\newcount\exnumber
\def\nextex{\bLP \the\sectionnumber.\the\exnumber \quad
 \advance\exnumber by 1}

\nopagenumbers
\headline={\ifnum\pageno=1 \hfill \else\hss{\tenrm--\folio--}\hss \fi}

\def\assignnumber#1#2{%
	\ifundefined{#1}\relax\else\message{#1 already defined}\fi
	\expandafter\xdef\csname#1\endcsname
 {\if-\the\sectionnumber\else\the\sectionnumber.\fi\the#2}%
	}%

\hyphenation{pre-print}

\let\sPP=\smallbreak

\let\bPP=\bigbreak
\def\LP{\par\noindent}
\def\sLP{\smallbreak\noindent}

\def\bLP{\bigbreak\noindent}

\def\Sec#1 #2 {\vskip0pt plus.1\vsize\penalty-250\vskip0pt plus-.1\vsize
	\bigbreak\bigskip
	\sectionnumber{#1} \equationnumber0\thnumber0
	\noindent{\bf\if-#1\else #1. \fi#2}\par
	\nobreak\smallskip\noindent}
\def\Subsec#1 #2 {\bLP{\bf
 	\if-\the\sectionnumber\else\the\sectionnumber.\fi#1. #2.}\quad}
\def\Proof{\goodbreak\sLP{\bf Proof}\quad}
\def\Ref{\vskip0pt plus.1\vsize\penalty-250\vskip0pt plus-.1\vsize
	\bigbreak\bigskip\leftline{\bf References}\nobreak\smallskip
	\frenchspacing}

\def\wrt{with respect to }

\def\halmos{\hbox{\vrule height0.31cm width0.01cm\vbox{\hrule height
 0.01cm width0.3cm \vskip0.29cm \hrule height 0.01cm width0.3cm}\vrule
 height0.31cm width 0.01cm}}
\def\hhalmos{{\unskip\nobreak\hfil\penalty50
	\quad\vadjust{}\nobreak\hfil\halmos
	\parfillskip=0pt\finalhyphendemerits=0\par}}
\def\nologo{\expandafter\let\csname logo\string @\endcsname=\empty}
\def\:{:\allowbreak }

\def\al{\alpha}
\def\be{\beta}
\def\ga{\gamma}
\def\de{\delta}
\def\ep{\varepsilon}

\def\th{\theta}

\def\la{\lambda}

\def\Ga{\Gamma}

\def\CC{{\Bbb C}}

\def\NN{{\Bbb N}}

\def\RR{{\Bbb R}}

\def\ZZ{{\Bbb Z}}

\def\FSC{{\cal C}}

\def\iy{\infty}

\def\const{{\rm const.}\,}

\def\Re{{\rm Re}\,}

\def\Zplus{\ZZ_+}

\def\LHS{left hand side}
\def\RHS{right hand side}

\def\rep{representation}

\let\wt=\widetilde

\def\eq#1{\relax
	\global\advance\equationnumber by 1
	\assignnumber{EN#1}\equationnumber
	{\rm (\csname EN#1\endcsname)}}
\def\eqtag#1{\relax\ifundefined{EN#1}\message{EN#1 undefined}{\sl (#1)}%
	\else(\csname EN#1\endcsname)\fi%
	}

\def\thname#1{\relax
	\global\advance\thnumber by 1
	\assignnumber{TH#1}\thnumber
	\csname TH#1\endcsname}
\def\Assumption#1 {\bLP{\bf Assumption \thname{#1}}\quad}
\def\Cor#1 {\bLP{\bf Corollary \thname{#1}}\quad}
\def\Def#1 {\bLP{\bf Definition \thname{#1}}\quad}
\def\Example#1 {\bLP{\bf Example \thname{#1}}\quad}
\def\Lemma#1 {\bLP{\bf Lemma \thname{#1}}\quad}
\def\Prop#1 {\bLP{\bf Proposition \thname{#1}}\quad}
\def\Remark#1 {\bLP{\bf Remark \thname{#1}}\quad}
\def\Theor#1 {\bLP{\bf Theorem \thname{#1}}\quad}
\def\thtag#1{\relax\ifundefined{TH#1}\message{TH#1 undefined}{\sl #1}%
	\else\csname TH#1\endcsname\fi}

\def\refitem#1 #2\par{\ifundefined{REF#1}
 \global\advance\refnum by1%
 \expandafter\xdef\csname REF#1\endcsname{\the\refnum}%
 \else\item{\ref{#1}}#2\sLP\fi}
\def\ref#1{\ifundefined{REF#1}\message{REF#1 is undefined}\else
	[\csname REF#1\endcsname]\fi}
\def\reff#1#2{\ifundefined{REF#1}\message{REF#1 is undefined}\else
	[\csname REF#1\endcsname, #2]\fi}
\refnum 0

\def\Pt{\wt P}

\font\titlefont=cmr10 scaled \magstep3
\font\rmten=cmr10 at 10truept
\font\ttten=cmtt10 at 10truept
\font\slten=cmsl10 at 10truept
\font\bften=cmbx10 at 10truept

\centerline{\titlefont
q-Special functions, a tutorial}
\bigskip
\bigskip
\centerline{TOM H. KOORNWINDER}
\bigskip
\midinsert
\narrower
{\rmten \noindent{\bften Abstract}\quad
A tutorial introduction is given to q-special functions and to
q-analogues of the classical orthogonal polynomials, up to the level
of Askey-Wilson polynomials.}
\endinsert

\refitem{AlCa65}

\refitem{And86}

\refitem{AnAs77}

\refitem{AnAs85}

\refitem{AsIs83}

\refitem{AsIs84}

\refitem{AsWi79}

\refitem{AsWi85}

\refitem{AtSu85}

\refitem{AtRaSu93}

\refitem{Bai35}

\refitem{Chi78}

\refitem{Erd1}

\refitem{Erd2}

\refitem{GaRa90}

\refitem{Hah49}

\refitem{Ism77}

\refitem{KaMi89}

\refitem{KoeSwa}

\refitem{Koo88}

\refitem{Koo89}

\refitem{Koo90}

\refitem{Koo93}

\refitem{KoSw}

\refitem{Lab90}

\refitem{Leo82}

\refitem{Mil89}

\refitem{Mim89}

\refitem{Moa81}

\refitem{Olv74}

\refitem{RaVe}

\refitem{Tit39}

\refitem{Wil80}

\refitem{Groe}

\par

\Sec0 {Introduction} \footnote{}
{\rmten
Korteweg-de Vries Institute, University of Amsterdam,
P.O. Box 94248, 1090 GE Amsterdam, The Netherlands;
{\slten email:} {\ttten T.H.Koornwinder@uva.nl}.
This paper is a slightly corrected version
(see in particular section 2.5) of sections 3 and~4 in
T. H. Koornwinder,
{\slten Compact quantum groups and q-special functions},
in {\slten Representations of Lie groups and quantum groups},
V. Baldoni \& M. A. Picardello (eds.),
Pitman Research Notes in Mathematics Series 311,
Longman Scientific \& Technical, pp. 46--128,
1994.}

It is the purpose of this paper to give a tutorial introduction to
$q$-hypergeometric functions and to orthogonal polynomials expressible
in terms of such functions.
An earlier version of this paper was written for
an intensive course on special functions aimed at Dutch graduate students,
it was elaborated during seminar lectures at Delft University of Technology,
and later it was part of the lecture notes
of my course on ``Quantum groups and $q$-special functions''
at the European School of Group Theory 1993, Trento, Italy.

I now describe the various sections in some more detail.
Section 1 gives an introduction to $q$-hypergeometric functions.
The more elementary $q$-special functions like $q$-exponential and
$q$-binomial series are treated in a rather self-contained way,
but for the higher $q$-hypergeometric functions some identities are given
without proof. The reader is referred, for instance, to the encyclopedic
treatise by Gasper \& Rahman \ref{GaRa90}.
Hopefully, this section succeeds to give the reader some feeling for
the subject and some impression of general techniques and ideas.

Section 2 gives an overview of the classical orthogonal
polynomials, where ``classical'' now means ``up to the level of
Askey-Wilson polynomials'' \ref{AsWi85}.
The section starts with the ``very classical'' situation of Jacobi, Laguerre and
Hermite polynomials and next  discusses the Askey tableau of classical
orthogonal polynomials (still for $q=1)$. Then the example of big $q$-Jacobi
polynomials is worked out in detail, as a demonstration how the main formulas
in this area can be neatly derived.
The section continues with the $q$-Hahn tableau and then gives a
self-contained introduction to the Askey-Wilson polynomials.
Both sections conclude with some exercises.

\sLP
{\sl Notation}\quad
The notations $\NN:=\{1,2,\ldots\}$ and
$\Zplus:=\{0,1,2,\ldots\}$ will be used.

\sLP
{\sl Acknowledgement}\quad
I thank Ren\'e Swarttouw, Jasper Stokman and Doug Bowman for commenting on
an earlier version of this paper.

\Sec1 {Basic hypergeometric functions}
A standard reference for this section is the book by
Gasper and Rahman \ref{GaRa90}.
In particular, the present section is quite parallel to their
introductory Chapter 1.
See also the useful compendia of formulas in the Appendices to that book.
The foreword to \ref{GaRa90}
by R. Askey gives a succinct historical introduction to the subject.
For first reading on the subject I can also recommend
Andrews \ref{And86}.

\Subsec1 {Preliminaries}
We start with briefly recalling the definition of the general hypergeometric
series (see Erd\'elyi e.a.\ \reff{Erd1}{Ch.\ 4} or Bailey
\reff{Bai35}{Ch.\ 2}).
For $a\in\CC$ the {\sl shifted factorial} or {\sl Pochhammer symbol} is
defined by
$(a)_0:=1$ and
$$
(a)_k:=a\,(a+1)\ldots(a+k-1),\quad k=1,2,\ldots\;.
$$
The general {\sl hypergeometric series} is defined by
$$
{}_rF_s\left[{
a_1,\ldots,a_r\atop b_1,\ldots,b_s};z\right]
:=\sum_{k=0}^\infty{(a_1)_k\ldots(a_r)_k\over (b_1)_k\ldots(b_s)_k\,k!}\,z^k.
\eqno\eq{116}
$$
Here $r,s\in \Zplus$ and the upper parameters $a_1,\ldots,a_r$,
the lower parameters $b_1,\ldots,b_s$ and the argument $z$ are in $\CC$.
However, in order to avoid zeros in the denominator, we require
that $b_1,\ldots,b_s\notin\{0,-1,-2,\ldots\}$.
If, for some $i=1,\ldots,r$, $a_i$ is a non-positive integer then the
series \eqtag{116} is terminating. Otherwise, we have an infinite
power series with radius of convergence equal to 0, 1 or $\iy$
according to whether $r-s-1>0$, $=0$ or $<0$, respectively.
On the \RHS\ of \eqtag{116} we have that
$$
{(k+1)\hbox{th term}\over
k\hbox{th term}}=
{(k+a_1)\,\ldots\,(k+a_r)\,z\over
(k+b_1)\,\ldots\,(k+b_s)\,(k+1)}
\eqno\eq{117}
$$
is rational in $k$.
Conversely, any rational function in $k$ can be written in the form of the
\RHS\ of \eqtag{117}.
Hence, any series
$\sum_{k=0}^\infty c_k$ with $c_0=1$ and $c_{k+1}/c_k$ rational in $k$ is of
the form of a hypergeometric series \eqtag{116}.

The cases ${}_0F_0$ and ${}_1F_0$ are elementary: exponential resp.\
binomial series. The case ${}_2F_1$ is the familiar Gaussian
hypergeometric series, cf.\ \reff{Erd1}{Ch.\ 2}.

We next give the definition of basic hypergeometric series.
For $a,q\in\CC$
define the {\sl $q$-shifted factorial\/} by
$(a;q)_0:=1$ and
$$
(a;q)_k:=(1-a)(1-aq)\ldots(1-aq^{k-1}),\quad k=1,2,\ldots\;.
$$
For $|q|<1$ put
$$
(a;q)_\infty:=\prod_{k=0}^\infty(1-aq^k).
$$
We also write
$$
(a_1,a_2,\ldots,a_r;q)_k:=(a_1;q)_k\,(a_2;q)_k\,\ldots\,(a_r;q)_k\,,\quad
k=0,1,2,\ldots\;\hbox{ or }\infty.
$$
Then a {\sl basic hypergeometric series\/} or {\sl $q$-hypergeometric series\/}
is defined by
$$
\displaylines{
\qquad{}_r\phi_s\left[{a_1,\ldots,a_r\atop b_1,\ldots,b_s};q,z\right]
={}_r\phi_s(a_1,\ldots,a_r;b_1,\ldots,b_s;q,z)\hfill
\cr
\hfill
:=\sum_{k=0}^\infty
{(a_1,\ldots,a_r;q)_k
\over(b_1,\ldots,b_s,q;q)_k}\,
\bigl((-1)^k\,q^{k(k-1)/2}\bigr)^{1+s-r}\,z^k,\quad
r,s\in\Zplus.
\qquad\eq{1.10}
\cr}
$$
On the \RHS\ of \eqtag{1.10} we have that
$$
{(k+1)\hbox{th term}\over
k\hbox{th term}}=
{(1-a_1q^k)\,\ldots\,(1-a_rq^k)\,(-q^k)^{1+s-r}\,z\over
(1-b_1q^k)\,\ldots\,(1-b_sq^k)\,(1-q^{k+1})}
\eqno\eq{1.20}
$$
is rational in $q^k$.
Conversely, any rational function in $q^k$ can be written in the form of the
\RHS\ of \eqtag{1.20}.
Hence, any series
$\sum_{k=0}^\infty c_k$ with $c_0=1$ and $c_{k+1}/c_k$ rational in $q^k$ is of
the form of a $q$-hypergeometric series \eqtag{1.10}.
This characterization is one explanation why we allow $q$ raised to a power
quadratic in $k$ in \eqtag{1.10}.

Because of the easily verified relation
$$
(a;q^{-1})_k=(-1)^k\,a^k\,q^{-k(k-1)/2}\,(a^{-1};q)_k\,,
$$
any series \eqtag{1.10} can be transformed into a series with base $q^{-1}$.
Hence, it is sufficient to study series \eqtag{1.10} with $|q|\le1$.
The tricky case $|q|=1$ has not yet been studied much and will be ignored by
us.
Therefore we may assume that $|q|<1$.
In fact, for convenience we will always assume that $0<q<1$, unless it is
otherwise stated.

In order to have a well-defined series \eqtag{1.10}, we require that
$$
b_1,\ldots,b_s\ne1,q^{-1},q^{-2},\ldots\;.
$$
The series \eqtag{1.10} will terminate iff, for some $i=1,\ldots,r$,
we have $a_i\in\{1,q^{-1},q^{-2},\ldots\}$.
If $a_i=q^{-n}$ ($n=0,1,2,\ldots$) then all terms in the series with
$k>n$ will vanish.
In the non-vanishing case, the convergence behaviour of \eqtag{1.10} can be
read off from \eqtag{1.20} by use of the ratio test.
We conclude:
$$
\hbox{convergence radius of \eqtag{1.10}}=
\cases{
\infty&if $r<s+1$,
\cr
1&if $r=s+1$,
\cr
0&if $r>s+1$.}
$$

We can view the $q$-shifted factorial as a $q$-analogue of the shifted
factorial
by the limit formula
$$
\lim_{q\to1}{(q^a;q)_k\over(1-q)^k}=(a)_k:=a(a+1)\ldots(a+k-1).
$$
Hence ${}_r\phi_s$ is a $q$-analogue of ${}_r F_s$ by the formal (termwise)
limit
$$
\lim_{q\uparrow1}{}_r\phi_s\left[{
q^{a_1},\ldots,q^{a_r}\atop q^{b_1},\ldots,q^{b_s}};q,(q-1)^{1+s-r}z\right]
={}_rF_s\left[{
a_1,\ldots,a_r\atop b_1,\ldots,b_s};z\right].
\eqno\eq{1.30}
$$
However, we get many other $q$-analogues of ${}_rF_s$ by adding upper or
lower parameters to the \LHS\ of \eqtag{1.30} which are equal to 0 or which
depend on $q$ in such a way that they tend to a limit $\ne1$ as
$q\uparrow1$.
Note that the notation \eqtag{1.10} has the drawback that some rescaling
is necessary before we can take limits for $q\to 1$.
On the other hand, parameters can be put equal to zero without problem in
\eqtag{1.10}, which would not be the case if we worked with
$a_1,\ldots,a_r$, $b_1,\ldots,b_s$ as in \eqtag{1.30}.
In general, $q$-hypergeometric series can be studied in their own right,
without much consideration for the $q=1$ limit case.
This philosophy is often (but not always) reflected in the notation
generally in use.

It is well known that the confluent ${}_1F_1$ hypergeometric
function can be obtained from the Gaussian ${}_2F_1$ hypergeometric
function by a limit process called {\sl confluence}.
A similar phenomenon occurs for $q$-series. Formally, by taking termwise
limits we have
$$
\lim_{a_r\to\infty} {}_r\phi_s\left[
{a_1,\ldots,a_r\atop b_1,\ldots,b_s};q,{z\over a_r}\right]
=
{}_{r-1}\phi_s\left[
{a_1,\ldots,a_{r-1}\atop b_1,\ldots,b_s};q,z\right].
\eqno\eq{1.25}
$$
In particular, this explains the particular choice of the
factor $\bigl((-1)^k\,q^{k(k-1)/2}\bigr)^{1+s-r}$ in \eqtag{1.10}.
If $r=s+1$ this factor is lacking and for $r<s+1$ it is naturally obtained
by the confluence process.
The proof of \eqtag{1.25} is by the following lemma.

\Lemma{1.5}
Let the complex numbers $a_k$, $k\in\Zplus$, satisfy the estimate
$|a_k|\le R^{-k}$ for some $R>0$. Let
$$
F(b;q,z):=\sum_{k=0}^\infty a_k\,(b;q)_k\,z^k,\quad |z|<R,\quad b\in\CC,
\quad 0<q<1.
$$
Then
$$
\lim_{b\to\infty}
F(b;q,z/b)=\sum_{k=0}^\infty a_k\,(-1)^k\,q^{k(k-1)/2}\,z^k,
\eqno\eq{1.27}
$$
uniformly for $z$ in compact subsets of $\CC$.

\Proof
For the $k$th term of the series on the \LHS\ of \eqtag{1.27} we have
$$
\eqalignno{
c_k
:=&
a_k\,(b;q)_k\,(z/b)^k
\cr
=&
a_k\,(b^{-1}-1)\,(b^{-1}-q)\ldots(b^{-1}-q^{k-1})\,z^k.
\cr}
$$
Now let $|b|>1$ and let $N\in\NN$ be such that
$q^{N-1}>|b|^{-1}\ge q^N$.
If $k\le N$ then
$$
|c_k|\le
|a_k|\,q^{k(k-1)/2}\,(2|z|)^k.
$$
If $k>N$ then
$$
\eqalignno{
|(b^{-1}-1)(b^{-1}-q)\ldots(b^{-1}-q^{k-1})|
&\le
2^k\,q^{N(N-1)/2}\,|b|^{N-k}
\cr
&\le
2^k\,|b|^{(1-N)/2}\,|b|^{N-k}
\cr
&\le
2^k\,|b|^{-k/2},
\cr}
$$
so
$$
|c_k|\le
|a_k|\,(2\,|z|\,|b|^{-1/2})^k.
$$
Now fix $M>0$.
Then, for each $\ep>0$ we can find $K>0$ and $B>1$ such that
$$
\sum_{k=K}^\infty|c_k|<\ep\quad
\hbox{if $|z|<M$ and $|b|>B$.}
$$
Combination with the termwise limit result completes the proof.\hhalmos

\bPP
It can be observed quite often in literature on $q$-special functions
that no rigorous limit proofs are given,
but only formal limit statements like \eqtag{1.25} and \eqtag{1.30}.
Sometimes, in heuristic reasonings, this is acceptable and productive for
quickly finding new results.
But in general, I would say that rigorous limit proofs have to be added.

Any terminating power series
$$
\sum_{k=0}^n c_k\,z^k
$$
can also be written as
$$
z^n\,\sum_{k=0}^n c_{n-k}\,(1/z)^k.
$$
When we want to do this for a terminating $q$-hypergeometric series,
we have to use that
$$
(a;q)_{n-k}=
{(a;q)_n\over(q^{n-1}a;q^{-1})_k}=
(a;q)_n\,{(-1)^k\,q^{k(k-1)/2}\,(a^{-1}q^{1-n})^k\over
(q^{1-n}a^{-1};q)_k}
$$
and
$$
{(q^{-n};q)_{n-k}\over(q;q)_{n-k}}=
{(q^{-n};q)_n\over(q;q)_n}\,{(q^n;q^{-1})_k\over(q^{-1};q^{-1})_k}
=(-1)^n\,q^{-n(n+1)/2}\,{(q^{-n};q)_k\over(q;q)_k}\,q^{(n+1)k}.
$$
Thus we obtain, for $n\in\Zplus$,
$$
\displaylines{
\quad{}_{s+1}\phi_s\left[{q^{-n},a_1,\ldots,a_s\atop b_1,\ldots,b_s};q,z\right]
=(-1)^n\,q^{-n(n+1)/2}\,{(a_1,\ldots,a_s;q)_n\over(b_1,\ldots,b_s;q)_n}\,z^n
\hfill
\cr
\hfill\times
{}_{s+1}\phi_s\left[
{q^{-n},q^{-n+1}b_1^{-1},\ldots,q^{-n+1}b_s^{-1}\atop
q^{-n+1}a_1^{-1},\ldots,q^{-n+1}a_s^{-1}};q,
{q^{n+1}b_1\ldots b_s\over a_1\ldots a_s z}\right].\quad
\eq{1.35}
\cr}
$$
Similar identities can be derived for other ${}_r\phi_s$ and for cases
where some of the parameters are 0.

Thus, any explicit evaluation of a terminating $q$-hypergeometric series
immediately implies a second one by inversion of the direction of summation
in the series, while any identity between two terminating $q$-hypergeometric
series implies three other identities.

\Subsec2 {The $q$-integral}
Standard operations of classical analysis like differentiation and
integration do not fit very well with $q$-hypergeometric series and can
be better replaced by $q$-derivative and $q$-integral.

The {\sl $q$-derivative} $D_q f$ of a function $f$ on an open real interval is
given by
$$
(D_qf)(x):={f(x)-f(qx)\over(1-q)x},\quad x\ne0,
\eqno\eq{1.36}
$$
and $(D_qf)(0):=f'(0)$ by continuity, provided $f'(0)$ exists.
Note that
$\lim_{q\uparrow1}(D_qf)(x)\allowbreak=f'(x)$ if $f$ is differentiable.
Note also that, analogous to $d/dx\;(1-x)^n=-n\,(1-x)^{n-1}$, we have
$$
f(x)=(x;q)_n\Longrightarrow
(D_qf)(x)=-{1-q^n\over1-q}\,(qx;q)_{n-1}.
\eqno\eq{1.37}
$$

Now recall that $0<q<1$.
If $D_qF=f$ and $f$ is continuous then, for real $a$,
$$
F(a)-F(0)=a\,(1-q)\,\sum_{k=0}^\infty f(aq^k)\,q^k.
\eqno\eq{1.38}
$$
This suggests the definition of the {\sl $q$-integral}
$$
\int_0^a f(x)\,d_qx:=a\,(1-q)\,\sum_{k=0}^\infty f(aq^k)\,q^k.
\eqno\eq{1.40}
$$
Note that it can be viewed as an infinite Riemann sum with nonequidistant
mesh widths.
In the limit, as $q\uparrow1$, the \RHS\ of \eqtag{1.40} will tend to the
classical integral
$\int_0^a f(x)\,dx$.

{}From \eqtag{1.38} we can also obtain $F(a)-F(b)$ expressed in terms of $f$.
This suggests the definition
$$
\int_a^b f(x)\,d_qx:=
\int_0^a f(x)\,d_qx-\int_0^b f(x)\,d_qx.
\eqno\eq{1.50}
$$
Note that \eqtag{1.40} and \eqtag{1.50} remain valid if $a$ or $b$ is
negative.

There is no unique canonical choice for the $q$-integral from 0 to $\infty$.
We will put
$$
\int_0^\infty f(x)\,d_qx:=
(1-q)\,\sum_{k=-\infty}^\infty f(q^k)\,q^k
$$
(provided the sum converges absolutely).
The other natural choices are then expressed by
$$a\,\int_0^\infty f(ax)\,d_qx=
a\,(1-q)\,\sum_{k=-\infty}^\infty f(aq^k)\,q^k,\quad a>0.
$$
Note that the above expression remains invariant when we replace $a$ by
$aq^n$ ($n\in\ZZ$).

As an example consider
$$
\int_0^1 x^\al\,d_qx=(1-q)\,\sum_{k=0}^\infty q^{k(\al+1)}=
{1-q\over1-q^{\al+1}}\,,\quad\Re\al>-1,
\eqno\eq{1.60}
$$
which tends, for $q\uparrow1$, to
$$
{1\over\al+1}=\int_0^1 x^\al\,dx.
\eqno\eq{1.70}
$$
{}From the point of view of explicitly computing definite integrals by
Riemann sum approximation this is much more efficient than the
approximation with equidistant mesh widths
$$
{1\over n}\,\sum_{k=1}^n\left({k\over n}\right)^\al.
$$
The $q$-approximation \eqtag{1.60} of the definite integral \eqtag{1.70}
(viewed as an area) goes essentially back to Fermat.

\Subsec3 {Elementary examples}
The {\sl $q$-binomial series} is defined by
$$
{}_1\phi_0(a;-;q,z):=\sum_{k=0}^\infty{(a;q)_k\over(q;q)_k}\,z^k,\quad
|z|<1.
\eqno\eq{1.80}
$$
Here the `$-$' in a ${}_r\phi_s$ expression denotes an empty parameter list.
The name ``$q$-binomial''
is justified since \eqtag{1.80}, with $a$ replaced by $q^a$,
formally tends, as $q\uparrow1$, to the
binomial series
$$
{}_1F_0(a;z):=\sum_{k=0}^\infty{(a)_k\over k!}\,z^k=(1-z)^{-a},\quad|z|<1.
\eqno\eq{1.90}
$$
The limit transition
is made rigorous in \reff{Koo90}{Appendix A}.
It becomes elementary in the case
of a terminating series ($a:=q^{-n}$ in \eqtag{1.80} and $a:=-n$ in
\eqtag{1.90}, where $n\in\Zplus$).

The $q$-analogue of the series evaluation in \eqtag{1.90} is as follows.

\Prop{1.10}
We have
$$
{}_1\phi_0(a;-;q,z)=
{(az;q)_\infty\over(z;q)_\infty}\,,\quad|z|<1.
\eqno\eq{1.100}
$$
In particular,
$$
{}_1\phi_0(q^{-n};-;q,z)=(q^{-n}z;q)_n\,.
\eqno\eq{1.102}
$$

\Proof
Put
$h_a(z):={}_1\phi_0(a;-;q,z)$.
Then
$$
(1-z)\,h_a(z)=(1-az)\,h_a(qz),\quad{\rm hence}\quad
h_a(z)={(1-az)\over(1-z)}\,h_a(qz).
$$
Iteration gives
$$
h_a(z)={(az;q)_n\over(z;q)_n}\,h_a(q^nz),\quad n\in\Zplus.
$$
Now use that $h_a$ is analytic and therefore continuous at 0 and use that
$h_a(0)=1$.
Thus, for $n\to\infty$ we obtain \eqtag{1.100}.\hhalmos

\bPP
The two most elementary $q$-hypergeometric series are the
two {\sl $q$-exponential series}
$$
\eqalignno{
e_q(z)&:={}_1\phi_0(0;-;q,z)=
\sum_{k=0}^\infty {z^k\over(q;q)_k}={1\over(z;q)_\infty},\quad|z|<1,
&\eq{1.110}
\cr
\noalign{\hbox{and}}
E_q(z)&:={}_0\phi_0(-;-;q,-z)=
\sum_{k=0}^\infty{q^{k(k-1)/2}\,z^k\over(q;q)_k}=(-z;q)_\infty,\quad z\in\CC.
&\eq{1.120}
\cr}
$$
The evaluation in \eqtag{1.110} is a specialization of
\eqtag{1.100}.
On the other hand, \eqtag{1.120} is a confluent limit of \eqtag{1.100}:
$$
\lim_{a\to\infty}{}_1\phi_0(a;-;q,-z/a)={}_0\phi_0(-;-;q,-z)
$$
(cf.\ \eqtag{1.30}).
The limit is uniform for $z$ in compact subsets of $\CC$, see
Lemma \thtag{1.5}.
Thus the evaluation in \eqtag{1.120} follows also from \eqtag{1.100}.

It follows from \eqtag{1.110} and \eqtag{1.120} that
$$
e_q(z)\,E_q(-z)=1.
\eqno\eq{1.130}
$$
This identity is a $q$-analogue of $e^z\,e^{-z}=1$. Indeed,
the two $q$-exponential series are $q$-analogues of the exponential
series by the limit formulas
$$
\lim_{q\uparrow1}E_q((1-q)z)=e^z=
\lim_{q\uparrow1}e_q((1-q)z).
$$
The first limit is uniform on compacta of $\CC$ by the majorization
$$
\left|
{q^{k(k-1)/2}\,(1-q)^k\,z^k\over(q;q)_k}\right|\le{|z|^k\over k!}\,.
$$
The second limit then follows by use of \eqtag{1.130}.

Although we assumed the convention $0<q<1$, it is often useful to find out
for a given $q$-hypergeometric series what will be obtained by
changing $q$ into $q^{-1}$ and then rewriting things again in base $q$.
This will establish a kind of duality for $q$-hypergeometric series.
For instance, we have
$$
e_{q^{-1}}(z)=E_q(-qz),
$$
which can be seen from the power series definitions.

The following four identities, including \eqtag{1.110}, are obtained from
each other by trivial rewriting.
$$
\eqalignno{
\sum_{k=0}^\infty{(1-q)^k\,z^k\over(q;q)_k}
&={1\over((1-q)z;q)_\infty}\,,
\cr
\sum_{k=0}^\infty{z^k\over(q;q)_k}&={1\over(z;q)_\infty}\,,
\cr
(1-q)^{1-b}\,\sum_{k=0}^\infty q^{kb}(q^{k+1};q)_\infty
&={(q;q)_\infty\over(q^b;q)_\infty\,(1-q)^{b-1}}\,,
\cr
\int_0^{(1-q)^{-1}}t^{b-1}\,((1-q)qt;q)_\infty\,d_qt
&={(q;q)_\infty\over(q^b;q)_\infty\,(1-q)^{b-1}}\,,\quad\Re b>0.
&\eq{1.140}
\cr}
$$
As $q\uparrow1$, the first identity tends to
$$
\sum_{k=0}^\infty{z^k\over k!}=e^z,
$$
while the \LHS\ of \eqtag{1.140} tends formally to
$$
\int_0^\infty t^{b-1}\,e^{-t}\,dt.
$$
Since this last integral can be evaluated as $\Ga(b)$, it is tempting
to consider the \RHS\ of \eqtag{1.140} as a the {\sl $q$-gamma function}.
Thus we put
$$
\Ga_q(z):=
{(q;q)_\infty\over(q^z;q)_\infty\,(1-q)^{z-1}}
=\int_0^{(1-q)^{-1}}t^{z-1}\,E_q(-(1-q)qt)\,d_qt,
\quad\Re z>0.
$$
where the last identity follows from \eqtag{1.140}.
It was proved in \reff{Koo90}{Appendix B} that
$$
\lim_{q\uparrow1}\Ga_q(z)=\Ga(z),\quad z\ne0,-1,-2,\ldots\;.
\eqno\eq{1.145}
$$

We have just seen an example how an identity for $q$-hypergeometric series
can have two completely different limit cases as $q\uparrow1$.
Of course, this is achieved by different rescaling.
In particular, reconsideration of a power series as a $q$-integral
is often helpful and suggestive for obtaining distinct limits.

Regarding $\Ga_q$ it can yet be remarked that it satisfies the
functional equation
$$
\Ga_q(z+1)={1-q^z\over1-q}\,\Ga_q(z)
$$
and that
$$
\Ga_q(n+1)={(q;q)_n\over(1-q)^n}\,,\quad n\in\Zplus.
$$

Similarly to the chain of equivalent identities including \eqtag{1.110}
we have a chain including \eqtag{1.100}:
$$
\eqalignno{
\sum_{k=0}^\infty{(q^a;q)_k\,z^k\over(q;q)_k}&=
{(q^az;q)_\infty\over(z;q)_\infty}\,,
\cr
(1-q)\,\sum_{k=0}^\infty q^{kb}\,
{(q^{k+1};q)_\infty\over(q^{k+a};q)_\infty}&=
{(1-q)\,(q,q^{a+b};q)_\infty\over
(q^a,q^b;q)_\infty}\,,
\cr
\int_0^1 t^{b-1}\,{(qt;q)_\infty\over(q^at;q)_\infty}\,d_qt&=
{\Ga_q(a)\,\Ga_q(b)\over\Ga_q(a+b)}\,,\quad\Re b>0.
&\eq{1.150}
\cr}
$$
In the limit, as $q\uparrow1$, the first identity tends to
$$
\sum_{k=0}^\infty{(a)_k\,z^k\over k!}=(1-z)^{-a},
$$
while \eqtag{1.150} formally tends to
$$
\int_0^1 t^{b-1}\,(1-t)^{a-1}\,dt={\Ga(a)\,\Ga(b)\over\Ga(a+b)}\,.
\eqno\eq{1.152}
$$
Thus \eqtag{1.150} can be considered as a $q$-beta integral and we define
the {\sl $q$-beta function\/} by
$$
B_q(a,b):={\Ga_q(a)\,\Ga_q(b)\over\Ga_q(a+b)}=
{(1-q)\,(q,q^{a+b};q)_\infty\over
(q^a,q^b;q)_\infty}=
\int_0^1 t^{b-1}\,{(qt;q)_\infty\over(q^at;q)_\infty}\,d_qt.
\eqno\eq{1.153}
$$

\Subsec4 {Heine's ${}_2\phi_1$ series}
Euler's integral \rep\ for the ${}_2F_1$ hypergeometric function
$$
\eqalign{
{}_2F_1(a,b;c;z)=
{\Ga(c)\over\Ga(b)\Ga(c-b)}\,
\int_0^1 t^{b-1}\,(1-t)^{c-b-1}\,(1-tz)^{-a}\,dt,\qquad&
\cr
\Re c>\Re b>0,\; |\arg(1-z)|<\pi,&
\cr}
\eqno\eq{1.155}
$$
(cf.\ \reff{Erd1}{2.1(10)})
has the following $q$-analogue due to Heine.
$$
\eqalign{
{}_2\phi_1(q^q,q^b;q^c;q,z)=
{\Ga_q(c)\over\Ga_q(b)\Ga_q(c-b)}\,
\int_0^1 t^{b-1}\,
{(tq;q)_\infty\over(tq^{c-b};q)_\infty}\,
{(tzq^a;q)_\infty\over(tz;q)_\infty}\,d_qt,&
\cr
\Re b>0,\quad |z|<1.&
\cr}
\eqno\eq{1.160}
$$
Note that the \LHS\ and \RHS\ of \eqtag{1.160} tend formally
to the corresponding sides of \eqtag{1.155}.
The proof of \eqtag{1.160} is also analogous to the proof of \eqtag{1.155}.
Expand $(tzq^a;q)_\infty/(tz;q)_\infty$ as a power series in $tz$
by \eqtag{1.80}, interchange summation and $q$-integration, and
evaluate the resulting $q$-integrals by \eqtag{1.150}.

If we rewrite the $q$-integral in \eqtag{1.160} as a series according to the
definition \eqtag{1.40},
and if we replace $q^a,q^b,q^c$ by $a,b,c$ then we obtain the
following transformation formula:
$$
{}_2\phi_1(a,b;c;q,z)=
{(az;q)_\infty\over(z;q)_\infty}\,{(b;q)_\infty\over(c;q)_\infty}\,
{}_2\phi_1(c/b,z;az;q,b).
\eqno\eq{1.170}
$$
Although \eqtag{1.160} and \eqtag{1.170} are equivalent, they look quite
different. In fact, in its form \eqtag{1.170} the identity has no classical
analogue. We see a new phenomenon, not occurring for ${}_2F_1$, namely that
the independent variable $z$ of the \LHS\ mixes on the \RHS\ with the
parameters. So, rather than having a function of $z$ with parameters
$a,b,c$, we deal with a function of $a,b,c,z$ satisfying certain
symmetries.

Just as in the classical case (cf.\ \reff{Erd1}{2.1(14)}),
substitution of some special value of $z$ in the $q$-integral \rep\
\eqtag{1.160} reduces it to a $q$-beta integral which can be explicitly
evaluated.
We obtain
$${}_2\phi_1(a,b;c;q,c/(ab))=
{(c/a,c/b;q)_\infty\over(c,c/(ab);q)_\infty}\,,\quad|c/(ab)|<1,
\eqno\eq{1.175}
$$
where the more relaxed bounds on the parameters are obtained by
analytic continuation.
The terminating case of \eqtag{1.175} is
$$
{}_2\phi_1(q^{-n},b;c;q,cq^n/b)=
{(c/b;q)_n\over(c;q)_n}\,,\quad n\in\Zplus.
\eqno\eq{1.177}
$$

The two fundamental transformation formulas
$$
\eqalignno{
{}_2F_1(a,b;c;z)&=
(1-z)^{-a}\,{}_2F_1(a,c-b;c;z/(z-1))
&\eq{1.178}
\cr
&=(1-z)^{c-a-b}\,{}_2F_1(c-a,c-b;c;z)
\cr}
$$
(cf.\ \reff{Erd1}{2.1(22) and (23)})
have the following $q$-analogues.
$$
\eqalignno{
{}_2\phi_1(a,b;c;q,z)
&={(az;q)_\infty\over(z;q)_\infty}\,
{}_2\phi_2(a,c/b;c,az;q,bz)
&\eq{1.180}
\cr
&={(abz/c;q)_\infty\over(z;q)_\infty}\,
{}_2\phi_1(c/a,c/b;c;q,abz/c).
&\eq{1.190}
\cr}
$$
Formula \eqtag{1.190} can be proved either by threefold iteration of
\eqtag{1.170} or by twofold iteration of \eqtag{1.180}.
The proof of \eqtag{1.180} is more involved.
Write both sides as power series in $z$.
Then make both sides into double series by substituting for
$(b;q)_k/(c;q)_k$ on the \LHS\ a terminating ${}_2\phi_1$
(cf.\ \eqtag{1.177}) and by substituting for
$(aq^kz;q)_\infty/(z;q)_\infty$ on the \RHS\ a $q$-binomial series
(cf.\ \eqtag{1.80}). The result follows by some rearrangement of series.
See \reff{GaRa90}{\S1.5} for the details.

Observe the difference between \eqtag{1.178} and its $q$-analogue
\eqtag{1.180}.
The argument $z/(z-1)$ in \eqtag{1.178} no longer occurs in \eqtag{1.180}
as a rational function of $z$, but the $z$-variable is distributed over
the argument and one of the lower parameters of the ${}_2\phi_2$.
Also we do not stay within the realm of ${}_2\phi_1$ functions.

Equation \eqtag{1.35} for $s:=1$ becomes
$$
\eqalign{
{}_2\phi_1(q^{-n},b;c;q,z)=
q^{-n(n+1)/2}\,{(b;q)_n\over(c;q)_n}\,(-z)^n\,
{}_2\phi_1\left(
q^{-n},{q^{-n+1}\over c};{q^{-n+1}\over b};q,{q^{n+1}c\over bz}\right),&
\cr
n\in\Zplus.&
\cr}
\eqno\eq{1.200}
$$
We may apply \eqtag{1.200} to the preceding evaluation and transformation
formulas for ${}_2\phi_1$ in order to obtain new ones
in the terminating case.
{}From \eqtag{1.177} we obtain
$$
{}_2\phi_1(q^{-n},b;c;q,q)=
{(c/b;q)_n\,b^n\over(c;q)_n}\,,\quad n\in\Zplus.
\eqno\eq{1.202}
$$
By inversion of direction of summation on both sides of \eqtag{1.180} 
(with $a:=q^{-n}$) we obtain
$$
{}_2\phi_1(q^{-n},b;c;q,z)=
{(c/b;q)_n\over(c;q)_n}\,
{}_3\phi_2\left[{
q^{-n},b,bzq^{-n}/c\atop bq^{1-n}/c,0};q,q\right],\quad n\in\Zplus.
\eqno\eq{1.204}
$$

A terminating ${}_2\phi_1$ can also be transformed into
a terminating ${}_3\phi_2$ with one of the upper parameters zero
(result of Jackson):
$$
{}_2\phi_1(q^{-n},b;c;q,z)=
(q^{-n}bz/c;q)_n\,
{}_3\phi_2\left[
{q^{-n},c/b,0\atop c,cqb^{-1}z^{-1}};q,q\right]
\eqno\eq{1.206}
$$
This formula can be proved by applying \eqtag{1.102}
to the factor
$(cq^{k+1}b^{-1}z^{-1};q)_{n-k}$
occuring in the $k$th term of the \RHS.
Then interchange summation in the resulting double sum
and substitute \eqtag{1.202} for the inner sum.
See \reff{Koo89}{ p.101} for the details.

\Subsec5 {A three-term transformation formula}
Formula \eqtag{1.200} has the following generalization for non-terminating
${}_2\phi_1$:
$$
\eqalignno{
&{}_2\phi_1(a,b;c;q,z)
+{(a,q/c,c/b,bz/q,q^2/bz;q)_\infty
\over
(c/q,aq/c,q/b,bz/c,cq/(bz);q)_\infty}\,
{}_2\phi_1\left({aq\over c},{bq\over c};{q^2\over c};q,z\right)\qquad
\cr
&=
{(abz/c,q/c,aq/b,cq/(abz);q)_\infty
\over
(bz/c,q/b,aq/c,cq/(bz);q)_\infty}\,
{}_2\phi_1\left(a,{aq\over c};{aq\over b};q,{cq\over abz}\right),
\quad \left|cq\over ab\right|<|z|<1.\qquad
&\eq{1.210}
\cr}
$$
This identity is a $q$-analogue of \reff{Erd1}{2.1(17)}.
In the following proposition we rewrite \eqtag{1.210} in an equivalent form
and next sketch the elegant proof due to
\ref{Mim89}.

\Prop{1.15}
Suppose $\{q^n\mid n\in\Zplus\}$ is disjoint from
$\{a^{-1}q^{-n},b^{-1}q^{-n}\mid n\in\Zplus\}$. Then
$$
\displaylines{
\quad(z,qz^{-1};q)_\infty\,{}_2\phi_1(a,b;c;q,z)
=
(az,qa^{-1}z^{-1};q)_\infty\,
{(c/a,b;q)_\infty\over(c,b/a;q)_\infty}\hfill
\cr
\hfill\times
{}_2\phi_1(a,qa/c;qa/b;q,qc/(abz))+(a\longleftrightarrow b).\quad
\eq{1.212}
\cr}
$$

\Proof
Consider the function
$$
F(w):=
{(a,b,cw,q,qwz^{-1},w^{-1}z;q)_\infty
\over
(aw,bw,c,w^{-1};q)_\infty}\,{1\over w}\,.
$$
Its residue at $q^n$ ($n\in\Zplus$) equals the $n^{\rm th}$ term of
the series on the \LHS\ of \eqtag{1.212}.
The negative of its residue at $a^{-1}q^{-n}$ ($n\in\Zplus$)
equals the $n^{\rm th}$ term of the first series on the \RHS\ of \eqtag{1.212},
and the residue at $b^{-1}q^{-n}$ is similarly related to the second series
on the \RHS.
Let $\FSC$ be a positively oriented closed curve around 0 in $\CC$ which
separates the two sets mentioned in the Proposition.
Then $(2\pi i)^{-1}\int_\FSC F(w)\,dw$ can be expressed 
in two ways as an infinite
sum of residues: either by letting the contour shrink to $\{0\}$ or by
blowing it up to $\{\infty\}$.\hhalmos

\bPP
If we put $z:=q$ in \eqtag{1.210} and substitute
\eqtag{1.175} then we obtain a generalization of \eqtag{1.202} for
non-terminating ${}_2\phi_1$:
$$
{}_2\phi_1(a,b;c;q,q)
+{(a,b,q/c;q)_\infty
\over
(aq/c,bq/c,c/q;q)_\infty}\,
{}_2\phi_1\left({aq\over c},{bq\over c};{q^2\over c};q,q\right)
=
{(abq/c,q/c;q)_\infty
\over
(aq/c,bq/c;q)_\infty}\,.
\eqno\eq{1.220}
$$
By \eqtag{1.50} this identity \eqtag{1.220} can be equivalently written
in $q$-integral form:
$$
\int_{qc^{-1}}^1
{(ct,qt;q)_\infty\over(at,bt;q)_\infty}\,d_qt
=
{(1-q)\,(abq/c,q/c,c,q;q)_\infty
\over
(aq/c,bq/c,a,b;q)_\infty}\,.
\eqno\eq{1.230}
$$

Replace in \eqtag{1.230} $c$, $a$, $b$ by $-q^c$, $-q^a$, $q^b$, respectively,
and let $q\uparrow1$. Then we obtain formally
$$
\int_{-1}^1(1+t)^{a-c}\,(1-t)^{b-1}\,dt=
{2^{b+a-c}\,\Ga(a-c+1)\,\Ga(b)
\over\Ga(a+b-c+1)}\,.
\eqno\eq{1.240}
$$
Note that, although it is trivial to obtain \eqtag{1.240} from \eqtag{1.152}
by an affine transformation of the integration variable,
their $q$-analogues \eqtag{1.150} and \eqtag{1.230} are by no means trivially
equivalent.

\Subsec6 {Bilateral series}
Put $a:=q$ in \eqtag{1.210} and substitute \eqtag{1.80}.
Then we can combine the two series from 0 to $\infty$ into a series from
$-\infty$ to $\infty$ without ``discontinuity'' in the summand:
$$
\sum_{k=-\infty}^\infty
{(b;q)_k\over(c;q)_k}\,z^k=
{(q,c/b,bz,q/(bz);q)_\infty\over
(c,q/b,z,c/(bz);q)_\infty}\,,\quad
|c/b|<|z|<1.
\eqno\eq{1.250}
$$
Here we have extended the definition of $(b;q)_k$ by
$$
(b;q)_k:={(b;q)_\infty\over(bq^k;q)_\infty}\,,\quad k\in\ZZ.
$$
Formula \eqtag{1.250} was first obtained by Ramanujan
({\sl Ramanujan's ${}_1\psi_1$-summation formula}).
It reduces to the $q$-binomial formula \eqtag{1.100} for $c:=q^n$
($n\in\Zplus$).
In fact, this observation can be used for a proof of \eqtag{1.250},
cf. \ref{Ism77}, \reff{And86}{Appendix C}.
Formula \eqtag{1.250} is a $q$-analogue of the explicit Fourier series
evaluation
$$
\eqalign{
\sum_{k=-\infty}^\infty{(b)_k\over(c)_k}\,e^{ik\th}=
{\Ga(c)\,\Ga(1-b)\over\Ga(c-b)}\,
e^{i(1-c)(\th-\pi)}\,(1-e^{i\th})^{c-b-1},\qquad&
\cr
0<\th<2\pi,\;\Re(c-b-1)>0,&
\cr}
$$
where $(b)_k:=\Ga(b+k)/\Ga(b)$,
also for $k=-1,-2,\ldots\;$.

Define {\sl bilateral $q$-hypergeometric series\/} by
$$
\eqalignno{
{}_r\psi_s&
\left[{a_1,\ldots,a_r\atop b_1,\ldots,b_s};q,z\right]\quad
\cr
:=&\sum_{k=-\infty}^\infty
{(a_1,\ldots,a_r;q)_k\over(b_1,\ldots,b_s;q)_k}\,
\bigl((-1)^k\,q^{k(k-1)/2}\bigr)^{s-r}\,z^k\quad
&\eq{1.260}
\cr
=&{}_{r+1}\phi_s\left[{
a_1,\ldots,a_r,q\atop b_1,\ldots,b_s};q,z\right]\quad
\cr
+&
{(b_1-q)\ldots(b_s-q)\over(a_1-q)\ldots(a_r-q)\,z}\,
{}_{s+1}\phi_s\left[{
q^2/b_1,\ldots,q^2/b_s,q\atop
q^2/a_1,\ldots,q^2/a_r,0,\ldots,0};q,
{b_1\ldots b_s\over a_1\ldots a_rz}\right],\quad
&\eq{1.270}
\cr}
$$
where $a_1,\ldots a_r,b_1,\ldots,b_s\ne0$ and $s\ge r$.
The Laurent series in \eqtag{1.260} is convergent for
$$
\eqalignno{
\left|{b_1\ldots b_s\over a_1\ldots a_r}\right|<|z|\quad
&{\rm if}\quad s>r
\cr
\noalign{\hbox{and for}}
\left|{b_1\ldots b_s\over a_1\ldots a_r}\right|<|z|<1\quad
&{\rm if}\quad s=r.
\cr}
$$
If some of the lower parameters in the ${}_r\psi_s$ are 0 then a suitable
confluent limit has to be taken in the ${}_{s+1}\phi_s$ of \eqtag{1.270}.

Thus we can write \eqtag{1.250} as
$$
{}_1\psi_1(b;c;q,z)=
{(q,c/b,bz,q/(bz);q)_\infty\over
(c,q/b,z,c/(bz);q)_\infty}\,,\quad
|c/b|<|z|<1.
\eqno\eq{1.280}
$$
Replace in \eqtag{1.280} $z$ by $z/b$,
substitute \eqtag{1.270},
let $b\to\infty$, and apply \eqtag{1.25}
Then we obtain
$$
{}_0\psi_1(-;c;q,z)
:=
\sum_{k=-\infty}^\infty
{(-1)^k\,q^{k(k-1)/2}\,z^k\over(c;q)_k}
=
{(q,z,q/z;q)_\infty\over(c,c/z;q)_\infty}\,,\quad|z|>|c|.
\eqno\eq{1.290}
$$
In particular, for $c=0$, we get the {\sl Jacobi triple product identity}
$$
{}_0\psi_1(-;0;q,z)
:=
\sum_{k=-\infty}^\infty(-1)^k\,q^{k(k-1)/2}\,z^k
=
(q,z,q/z;q)_\infty\,,\quad z\ne0.
\eqno\eq{1.300}
$$

The series in \eqtag{1.300} is essentially a {\sl theta function}.
With the notation \reff{Erd2}{13.19(9)} we get
$$
\eqalignno{
\th_4(x;q)
:=&
\sum_{k=-\infty}^\infty(-1)^k\,q^{k^2}\,e^{2\pi ikx}
\cr
=&
{}_0\psi_1(-;0;q^2,q\,e^{2\pi ix})
\cr
=&
(q^2,q\,e^{2\pi ix},q\,e^{-2\pi ix};q^2)_\infty
\cr
=&
\prod_{k=1}^\infty(1-q^{2k})\,
(1-2q^{2k-1}\,\cos(2\pi x)+q^{4k-2}),
\cr}
$$
and similarly for the other theta functions $\th_i(x;q)$ ($i=1,2,3$).

\Subsec7 {The $q$-hypergeometric $q$-difference equation}
Just as the
{\sl hypergeometric differential equation}
$$
z\,(1-z)\,u''(z)+(c-(a+b+1)z)\,u'(z)-a\,b\,u(z)=0
$$
(cf.\ \reff{Erd1}{2.1(1)}) has particular solutions
$$
u_1(z):={}_2F_1(a,b;c;z),\quad
u_2(z):=z^{1-c}\,{}_2F_1(a-c+1,b-c+1;2-c;z),
$$
the {\sl $q$-hypergeometric $q$-difference equation\/}
$$
\displaylines{
\quad z\,(q^c-q^{a+b+1}z)\,(D_q^2u)(z)+
\left[{1-q^c\over1-q}-
\left(q^b\,{1-q^a\over1-q}+q^a\,{1-q^{b+1}\over1-q}\right)\,z\right]\,
(D_qu)(z)\hfill
\cr
\hfill-\,{1-q^a\over1-q}\,{1-q^b\over1-q}\,u(z)=0\quad
\eq{1.310}
\cr}
$$
has particular solutions
$$
\eqalignno{
u_1(z)&:={}_2\phi_1(q^a,q^b;q^c;q,z),
&\eq{1.312}
\cr
u_2(z)&:=z^{1-c}\,{}_2\phi_1(q^{1+a-c},q^{1+b-c};q^{2-c};q,z).
&\eq{1.314}
\cr}
$$
There is an underlying theory of $q$-difference equations with regular
singularities, similarly to the theory of differential equations with regular
singularities discussed for instance in Olver \reff{Olv74}{Ch. 5}.
It is not difficult to prove the following proposition.

\Prop{1.20}
Let $A(z):=\sum_{k=0}^\infty a_k\,z^k$ and
$B(z):=\sum_{k=0}^\infty b_k\,z^k$ be convergent power series.
Let $\la\in\CC$ be such that
$$
{(1-q^{\la+k})\,(1-q^{\la+k-1})\over(1-q)^2}+a_0\,{1-q^{\la+k}\over1-q}+b_0
\quad\cases{
=0,&$k=0$,
\cr
\ne0,&$k=1,2,\ldots\;$.
\cr}
\eqno\eq{1.316}
$$
Then the $q$-difference equation
$$
z^2\,(D_q^2u)(z)+z\,A(z)\,(D_qu)(z)+B(z)\,u(z)=0
\eqno\eq{1.320}
$$
has an (up to a constant factor) unique solution of the form
$$
u(z)=\sum_{k=0}^\infty c_k\,z^{\la+k}.
\eqno\eq{1.330}
$$

\bPP
Note that \eqtag{1.310} can be rewritten in the form \eqtag{1.320} with
$a_0=(q^{-c}-1)/(1-q)$, $b_0=0$,
so \eqtag{1.316} has solutions $\la=0$ and $-c+1$
(mod $(2\pi i\,\log q^{-1})\ZZ$)
provided $c\notin\ZZ$ (mod $(2\pi i\,\log q^{-1})\ZZ$).
For the coefficients $c_k$ in \eqtag{1.330} we find the recursion
$$
{c_{k+1}\over c_k}=
{(1-q^{a+\la+k})\,(1-q^{b+\la+k})\over
(1-q^{c+\la+k})\,(1-q^{\la+k+1})}\,.
$$
Thus we obtain solutions $u_1$, $u_2$ as given in \eqtag{1.312},
\eqtag{1.314}.

Proposition \thtag{1.20} can also be applied to the case $z=\infty$
of \eqtag{1.310}.
Just make the transformation $z\mapsto z^{-1}$ in \eqtag{1.310}.
One solution then obtained is
$$
u_3(z):=
z^{-a}\,{}_2\phi_1\left[{
q^a,q^{a-c+1}\atop q^{a-b+1}};q,q^{-a-b+c+1}\,z^{-1}\right].
$$
Now \eqtag{1.210} can be rewritten as
$$
\displaylines{
\quad u_1(z)+
{(q^a,q^{1-c},q^{c-b};q)_\infty\over
(q^{c-1},q^{a-c+1},q^{1-b};q)_\infty}
\,
{(q^{b-1}z,q^{2-b}z^{-1};q)_\infty\,z^{c-1}\over
(q^{b-c}z,q^{c-b+1}z^{-1};q)_\infty}\,u_2(z)\hfill
\cr
\hfill
={(q^{1-c},q^{a-b+1};q)_\infty\over
(q^{1-b},q^{a-c+1};q)_\infty}\,
{(q^{a+b-c}z,q^{c-a-b+1}z^{-1};q)_\infty\,z^a\over
(q^{b-c}z,q^{c-b+1}z^{-1};q)_\infty}\,u_3(z).\quad\eq{1.340}
\cr}
$$
Note that $u_3$ is not a linear combination of $u_1$ and $u_2$, as the
coefficients of $u_2$ and $u_3$ in \eqtag{1.340} depend on $z$.
However, since an expression of the form
$$
z\mapsto
{(q^\al z,q^{1-\al}z^{-1};q)_\infty\,z^{\al-\be}\over
(q^\be z,q^{1-\be}z^{-1};q)_\infty}
$$
is invariant under transformations $z\mapsto qz$,
each term in \eqtag{1.340} is a solution of \eqtag{1.310}.
Thus everything works fine when we restrict ourselves to a subset of the form
$\{z_0\,q^k\mid k\in\ZZ\}$.

The $q$-analogue of the regular singularity at $z=1$ for the ordinary
hypergeometric differential equation has to be treated in a different way.
We can rewrite \eqtag{1.310} as
$$
(q^c-q^{a+b}z)\,u(qz)+(-q^c-q+(q^a+q^b)z)\,u(z)+(q-z)\,u(q^{-1}z)=0.
$$
It can be expected that the points $z=q^{c-a-b}$ and $z=q$,
where the coefficient of $u(qz)$ respectively $u(q^{-1}z)$ vanishes,
will replace the classical regular singularity $z=1$, see also
\eqtag{1.175}, \eqtag{1.202} and \eqtag{1.220}.
A systematic theory for such singularities has not yet been developed.

\Subsec8 {$q$-Bessel functions}
The classical {\sl Bessel function\/} is defined by
$$
\eqalignno{
J_\nu(z):=&(z/2)^\nu\,
\sum_{k=0}^\infty{(-1)^k\over\Ga(\nu+k+1)\,k!}\,(z/2)^{2k}
\cr
=&{(z/2)^\nu\over\Ga(\nu+1)}\,{}_0F_1(-;\nu+1;-z^2/4),
\cr}
$$
cf.\ \reff{Erd2}{Ch.\ 7}).
Jackson (1905) introduced two $q$-analogues of the Bessel function:
$$
\eqalignno{
J_\nu^{(1)}(z;q)
:=&
{(q^{\nu+1};q)_\infty\over(q;q)_\infty}\,(z/2)^\nu\,
{}_2\phi_1(0,0;q^{\nu+1};q,-z^2/4)
\cr
=&
(z/2)^\nu\,
\sum_{k=0}^\infty{(q^{\nu+k+1};q)_\infty\over(q;q)_\infty}\,
{(-1)^k\,(z/2)^{2k}\over(q;q)_k}\,,&\eq{114}
\cr
J_\nu^{(2)}(z;q):=&
{(q^{\nu+1};q)_\infty\over(q;q)_\infty}\,(z/2)^\nu\,
{}_0\phi_1(-;q^{\nu+1};q,-q^{\nu+1}z^2/4)
\cr
=&
(z/2)^\nu\,\sum_{k=0}^\infty{(q^{\nu+k+1};q)_\infty\over(q;q)_\infty}\,
{q^{k(k+\nu)}\,(-1)^k\,(z/2)^{2k}\over(q;q)_k}\,.&\eq{115}
\cr}
$$
Formally we have
$$
\lim_{q\uparrow1}J_\nu^{(i)}((1-q)z;q)=J_\nu(z),\quad i=1,2
$$
(cf.\ \eqtag{1.145}).
The two $q$-Bessel functions $J_\nu^{(i)}(z;q)$ can be simply expressed in terms
of each other.
{}From \eqtag{1.190} and \eqtag{1.25} we obtain
$$
{}_2\phi_1(0,b;c;q,z)=
{1\over(z;q)_\infty}\,{}_1\phi_1(c/b;c;q,bz).
$$
A further confluence with $b\to0$ yields
$$
{}_2\phi_1(0,0;c;q,z)=
{1\over(z;q)_\infty}\,{}_0\phi_1(-;c;q,cz).
$$
This can be rewritten as
$$
J_\nu^{(2)}(z;q)=(-z^2/4;q)_\infty\,J_\nu^{(1)}(z;q).
$$

Yet another $q$-analogue of the Bessel function is as follows.
$$
\eqalignno{
J_\nu(z;q)
:=&
{(q^{\nu+1};q)_\infty\over(q;q)_\infty}\,z^\nu\,
{}_1\phi_1(0;q^{\nu+1};q,qz^2)
\cr
=&
z^\nu\,\sum_{k=0}^\infty{(q^{\nu+k+1};q)_\infty\over(q;q)_\infty}\,
{(-1)^k\,q^{k(k+1)/2}\,z^{2k}\over(q;q)_k}\,.
&\eq{1.380}
\cr}
$$
Formally we have
$$
\lim_{q\uparrow1}J_\nu((1-q^{1/2})z;q)=J_\nu(z).
$$
This so-called {\sl Jackson's third $q$-Bessel function} is not simply
related to the
$q$-Bessel functions \eqtag{114}, \eqtag{115}, but they were all introduced
by Jackson.
Koornwinder \& Swarttouw \ref{KoSw} gave
a satisfactory $q$-analogue of the Hankel transform in terms
of the $q$-Bessel function \eqtag{1.380}.

\Subsec9 {Various results}
Goursat's list of
{\sl quadratic transformations\/} for Gaussian hypergeometric functions
can be found in \reff{Erd1}{\S2.11}.
In a recent paper Rahman \& Verma \ref{RaVe}
have given a full list of $q$-analogues of Goursat's table.
However, all their formulas involve on at least one of both sides
an ${}_8\phi_7$ series.
Moreover, a good foundation from the theory of
$q$-difference equations is not yet available.
For terminating series many of their ${}_8\phi_7$'s will simplify to
${}_4\phi_3$'s.
In section 2 we will meet some natural examples of these transformations
coming from orthogonal polynomials.

An important part of the book by Gasper \& Rahman \ref{GaRa90}
deals with the derivation of summation and transformation formulas
of ${}_{s+1}\phi_s$ functions with $s>1$.
A simple example is the {\sl $q$-Saalsch\"utz\/} formula
$$
{}_3\phi_2(a,b,q^{-n};c,abc^{-1}q^{1-n};q,q)
=
{(c/a,c/b;q)_n\over(c,c/(ab);q)_n}\,,\quad n\in\Zplus,
\eqno\eq{1.400}
$$
which follows easily from \eqtag{1.190}
by expanding the quotient on the \RHS\ of \eqtag{1.190} with the aid of
\eqtag{1.100}, and next comparing equal powers of $z$ at both sides of
\eqtag{1.190}.
Formula \eqtag{1.400} is the $q$-analogue of the
{\sl Pfaff-Saalsch\"utz formula}
$$
{}_3F_2(a,b,-n;c,1+a+b-c-n;1)=
{(c-a)_n\,(c-b)_n\over (c)_n\,(c-a-b)_n}\,,\quad n\in\Zplus,
$$
which can be proved in an analogous way as \eqtag{1.400}.

An example of a much more involved transformation formula, which has
many important special cases, is {\sl Watson's transformation formula}
$$
\displaylines{
\qquad{}_8\phi_7\left[
{a,qa^{1/2},-qa^{1/2},b,c,d,e,q^{-n}
\atop
a^{1/2},-a^{1/2},aq/b,aq/c,aq/d,aq/e,aq^{n+1}};q,{a^2q^{2+n}\over bcde}\right]
\hfill
\cr
\hfill
={(aq,aq/(de);q)_n\over(aq/d,aq/e;q)_n}\,
{}_4\phi_3\left[{
q^{-n},d,e,aq/(bc)
\atop
aq/b,aq/c,deq^{-n}/a};q,q\right],\quad n\in\Zplus\,,\qquad
\eq{1.410}
\cr}
$$
cf.\ Gasper \& Rahman \reff{GaRa90}{\S2.5}.
Then, for $a^2q^{n+1}=bcde$, the \RHS\ can be evaluated by use of
\eqtag{1.400}. The resulting evaluation
$$
{(aq,aq/(bc),aq/(bd),aq/(cd);q)_n
\over
(aq/b,aq/c,aq/d,aq/(bcd);q)_n}
$$
of the \LHS\ of \eqtag{1.410} subject to the given relation between
$a,b,c,d,e,q$ is called
{\sl Jackson's summation formula}, cf.\ \reff{GaRa90}{\S2.6}.

The famous {\sl Rogers-Ramanujan identities}
$$
\eqalignno{
&{}_0\phi_1(-;0;q,q)
:=\sum_{k=0}^\infty{q^{k^2}\over(q;q)_k}
=
{(q^2,q^3,q^5;q^5)_\infty
\over(q;q)_\infty}\,,
\cr
&{}_0\phi_1(-;0;q,q^2)
:=\sum_{k=0}^\infty{q^{k(k+1)}\over(q;q)_k}
=
{(q,q^4,q^5;q^5)_\infty\over(q;q)_\infty}
\cr}
$$
have been proved in many different ways
(cf.\ Andrews \ref{And86}), not only analytically but also
by an interpretation in combinatorics or in the framework of Kac-Moody
algebras.
A quick analytic proof starting from \eqtag{1.410} is described
in \reff{GaRa90}{\S2.7}.

\bPP\goodbreak\noindent \exnumber1
{\bf Exercises to \S\the\sectionnumber}
\par\nobreak\smallskip\noindent\the\sectionnumber.\the\exnumber\quad
\advance\exnumber by 1
Prove that
$$
(a;q)_n=(-a)^n\,q^{n(n-1)/2}\,(a^{-1}q^{1-n};q)_n.
$$

\nextex
Prove that
$$
\eqalignno{
(a;q)_{2n}&=(a;q^2)_n\,(aq;q^2)_n,
\cr
(a^2;q^2)_n&=(a;q)_n\,(-a;q)_n,
\cr
(a;q)_\infty&=(a^{1/2},-a^{1/2},(aq)^{1/2},-(aq)^{1/2};q)_\infty.
\cr}
$$

\nextex
Prove the following identity of Euler:
$$
(-q;q)_\infty\,(q;q^2)_\infty=1.
$$

\nextex
Prove that
$$
\sum_{l=-\infty}^\infty
{(-1)^{l-m}\,q^{(l-m)(l-m-1)/2}\over
\Ga_q(n-l+1)\,\Ga_q(l-m+1)}
=\de_{n,m},\quad 0<q\le1.
$$
Do it first for $q=1$. Start for instance with $e^z\,e^{-z}=1$ and,
in the $q$-case, with $e_q(z)\,E_q(-z)=1$.

\nextex
Let
$$
{1\over(q;q)_\infty}=\sum_{n=0}^\infty a_n\,q^n
$$
be the power series expansion of the \LHS\ in terms of $q$.
Show that $a_n$ equals the number of partitions of $n$.
\LP
Show also that the coefficient $b_{k,n}$ in
$$
{q^k\over(q;q)_k}=\sum_{n=k}^\infty b_{k,n}\,q^n
$$
is the number of partitions of $n$ with highest part $k$.
Give now a partition theoretic proof of the identity
$$
{1\over(q;q)_\infty}=\sum_{k=0}^\infty{q^k\over(q;q)_k}\,.
$$

\nextex
In the same way, give a partition theoretic proof of the identity
$$
{1\over(q^m;q)_\infty}=\sum_{k=0}^\infty{q^{mk}\over(q;q)_k}\,,\quad m\in\NN.
$$

\nextex
Give also a partition theoretic proof of
$$
\sum_{k=0}^\infty {q^{k(k-1)/2}q^{mk}\over(q;q)_k}=(-q^m;q)_\infty,\quad
m\in\NN.
$$
(Consider the problem first for $m=1$.)

\nextex
Let $GF(p)$ be the finite field with $p$ elements.
Let $A$ be a $n\times n$ matrix with entries chosen independently and
at random from $GF(p)$, with equal probability for the field elements to be
chosen. Let $q:=p^{-1}$.
Prove that the probability that $A$ is invertible is
$(q;q)_n$.
(See SIAM News 23 (1990) no.6, p.8.)

\nextex
Let $GF(p)$ be as in the previous exercise.
Let $V$ be an $n$-dimensional vector space over $GF(p)$.
Prove that the number of $k$-dimensional linear subspaces of $V$ equals
the $q$-binomial coefficient
$$
\left[n\atop k\right]_p:={(p;p)_n\over(p;p)_k\,(p;p)_{n-k}}.
$$

\nextex
Show that
$$
(ab;q)_n=\sum_{k=0}^n
\left[n\atop k\right]_q b^k\,(a;q)_k\,(b;q)_{n-k}.
$$
Show that both Newton's binomial formula and the formula
$$
(a+b)_n=\sum_{k=0}^n{n\choose k}(a)_k\,(b)_{n-k}
$$
are limit cases of the above formula.

\Sec2 {$q$-Analogues of the classical orthogonal polynomials}
Originally by classical orthogonal polynomials
were only meant the three families of Jacobi, Laguerre and Hermite
polynomials, but recent insights consider a much bigger class of polynomials
as ``classical''.
On the one hand, there is an extension to hypergeometric orthogonal
polynomials up to the ${}_4F_3$ level and including certain discrete
orthogonal polynomials.
These are brought together in the Askey tableau, cf.{} Table 1.
On the other hand there are $q$-analogues of all the families in the
Askey tableau,
often several $q$-analogues for one classical family
(cf.{} Table 2 for some of them).
The master class of all these $q$-analogues is formed by the celebrated
Askey-Wilson polynomials.
They contain all other families described in this chapter as special cases
or limit cases.
Good references for this chapter are Andrews \& Askey
\ref{AnAs85}
and Askey \& Wilson \ref{AsWi85}.
See Koekoek \& Swarttouw \ref{KoeSwa}
for a quite comprehensive list of formulas.
See also Atakishiyev, Rahman \& Suslov \ref{AtRaSu93} for a somewhat different
approach.

Some parts of this section contain surveys without many proofs.
However, the subsections
\the\sectionnumber.3 and \the\sectionnumber.4 on big and little
$q$-Jacobi polynomials and
\the\sectionnumber.5 and \the\sectionnumber.6 on the
Askey-Wilson integral and related polynomials are
rather self-contained.

\Subsec1 {Very classical orthogonal polynomials}
An introduction to the traditional classical orthogonal polynomials
can be found, for instance, in \reff{Erd2}{Ch. 10}.
One possible characterization is as systems of orthogonal polynomials
$\{p_n\}_{n=0,1,2,\ldots}$ which are eigenfunctions of a second order
differential operator not involving $n$
with eigenvalues $\la_n$ depending on $n$:
$$
a(x)\,p_n''(x)+b(x)\,p_n'(x)+c(x)\,p_n(x)=\la_n\,p_n(x).
\eqno\eq{2.5}
$$
Because we will extend the definition of classical orthogonal polynomials
in this section, we will call the orthogonal polynomials satisfying
\eqtag{2.5} {\sl very classical orthogonal polynomials.}
The classification shows that, up to an affine transformation of the
independent variable, the only cases are as follows.
\sLP
{\sl Jacobi polynomials}
$$
P_n^{(\al,\be)}(x):={(\al+1)_n\over n!}\,
{}_2F_1(-n,n+\al+\be+1;\al+1;(1-x)/2),\quad\al,\be>-1,
\eqno\eq{2.10}
$$
orthogonal on $[-1,1]$ \wrt\ the measure
$(1-x)^\al\,(1+x)^\be\,dx$;
\sLP
{\sl Laguerre polynomials}
$$
L_n^\al(x):={(\al+1)_n\over n!}\,{}_1F_1(-n;\al+1;x),\quad\al>-1,
\eqno\eq{2.20}
$$
orthogonal on $[0,\infty)$ \wrt\ the measure $x^\al\,e^{-x}\,dx$;
\sLP
{\sl Hermite polynomials}
$$
H_n(x):=(2x)^n\,{}_2F_0(-n/2,(1-n)/2;-;-x^{-2}),
\eqno\eq{2.30}
$$
orthogonal on $(-\infty,\infty)$ \wrt\ the measure $e^{-x^2}\,dx$.

In fact, Jacobi polynomials are the generic case here,
while the other two classes are limit cases of Jacobi polynomials:
$$
\eqalignno{
L_n^\al(x)&=\lim_{\be\to\infty}P_n^{(\al,\be)}(1-2x/\be),
&\eq{2.40}
\cr
H_n(x)&=
2^n\,n!\,\lim_{\al\to\infty}\al^{-n/2}\,P_n^{(\al,\al)}(\al^{-1/2}x).
&\eq{2.50}
\cr}
$$
Hermite polynomials are also limit cases of Laguerre
polynomials:
$$
H_n(x)=(-1)^n\,2^{n/2}\,n!\,\lim_{\al\to\infty}
\al^{-n/2}\,L_n^\al((2\al)^{1/2}x+\al).
\eqno\eq{2.60}
$$
The limit \eqtag{2.40} is immediate from \eqtag{2.10} and \eqtag{2.20}.
The limit \eqtag{2.50} follows from \eqtag{2.30} and a similar series
representation for {\sl Gegenbauer\/} or {\sl ultraspherical
polynomials\/} (special Jacobi polynomials with $\al=\be$):
$$
P_n^{(\al,\al)}(x)=
{(\al+1)_n\,(\al+1/2)_n\over(2\al+1)_n\,n!}\,(2x)^n\,
{}_2F_1(-n/2,(1-n)/2;-\al-n+1/2;x^{-2}),
$$
cf.\ \reff{Erd2}{10.9(4) and (18)}.
The limit \eqtag{2.60} cannot be easily derived by comparison of series
representations.
One method of proof is to rewrite the three term recurrence relation
for Laguerre polynomials
$$
(n+1)\,L_{n+1}^\al(x)-(2n+\al+1-x)\,L_n^\al(x)+(n+\al)\,L_{n-1}^\al(x)=0
$$
(cf.\ \reff{Erd2}{10.12(8)})
in terms of the polynomials in $x$ given by the right hand side of
\eqtag{2.60} and, next, to compare it with the three term recurrence relation
for Hermite polynomials (cf.\ \reff{Erd2}{10.13(10)})
$$
H_{n+1}(x)-2x\,H_n(x)+2n\,H_{n-1}(x)=0.
$$

Very classical orthogonal polynomials have other characterizations,
for instance by the existence of a Rodrigues type formula or by the fact
that the first derivatives again form a system of orthogonal polynomials
(cf.\ \reff{Erd2}{\S10.6}).
Here we want to point out that, associated with the last two characterizations,
there is a pair of differential recurrence relations from
which many of the basic properties of the polynomials can be easily
derived.
For instance, for Jacobi polynomials we have the pair
$$
\displaylines{
\hfill{d\over dx}\,P_n^{(\al,\be)}(x)=
{n+\al+\be+1\over 2}\,P_{n-1}^{(\al+1,\be+1)}(x),\hfill\eq{2.70}
\cr
\hfill(1-x)^{-\al}\,(1+x)^{-\be}\,
{d\over dx}\,\bigl((1-x)^{\al+1}\,(1+x)^{\be+1}\,
P_{n-1}^{(\al+1,\be+1)}(x)\bigr)
=
-2n\,P_n^{(\al,\be)}(x).
\hfill\eq{2.80}
\cr}
$$
The differential operators in \eqtag{2.70}, \eqtag{2.80} are called
{\sl shift operators\/} because of their parameter shifting property.
The differential operator in \eqtag{2.70} followed by the one in \eqtag{2.80}
yields the second order differential operator of which the Jacobi polynomials
are eigenfunctions.
If we would have defined Jacobi polynomials only by their orthogonality
property, not by their explicit expression, then we would have already been
able to derive \eqtag{2.70}, \eqtag{2.80} up to constant factors just by the
remarks that the operators $D_-$ in \eqtag{2.70} and $D_+^{(\al,\be)}$ in
\eqtag{2.80} satisfy
$$
\displaylines{
\qquad\int_{-1}^1(D_-f)(x)\,g(x)\,(1-x)^{\al+1}\,(1+x)^{\be+1}\,dx\hfill
\cr
\hfill=-\int_{-1}^1f(x)\,(D_+^{(\al,\be)}g)(x)\,(1-x)^\al\,(1+x)^\be\,dx\qquad
\eq{2.85}
\cr}
$$
and that $D_-$ sends polynomials of degree $n$ to polynomials of degree $n-1$,
while $D_+^{(\al,\be)}$ sends polynomials of degree $n-1$ to polynomials of
degree $n$.

The same idea can be applied again and again for the more general
classical orthogonal polynomials we will discuss in this section.

It follows from \eqtag{2.70}, \eqtag{2.80} and \eqtag{2.85}
that
$$
\int_{-1}^1\bigl(P_n^{(\al,\be)}(x)\bigr)^2\,
(1-x)^\al\,(1+x)^\be\,dx=
\const\int_{-1}^1(1-x)^{\al+n}\,(1+x)^{\be+n}\,dx,
$$
where the constant on the \RHS\ can be easily computed.
What is left for computation is a beta integral, which of course is elementary.
However, we emphasize this reduction of computation of quadratic norms of
classical orthogonal polynomials to computation of the integral of the weight
function with shifted parameter, because this phenomenon will also return
in the generalizations.
Usually, the computation of the integral of the weight function is the only
nontrivial part.
On the other hand, if we have some deep evaluation of a definite integral
with positive integrand, then it is worth to explore the polynomials
being orthogonal \wrt the weight function given by the integrand.

\Subsec2 {The Askey tableau}
Similar to the classification discussed in \S2.1,
one can classify
all systems of orthogonal polynomials $\{p_n\}$ which are
eigenfunctions of a second order difference operator:
$$
a(x)\,p_n(x+1)+b(x)\,p_n(x)+c(x)\,p_n(x-1)=\la_n\,p_n(x).
\eqno\eq{2.90}
$$
Here the definition of orthogonal polynomials is relaxed somewhat.
We include the possibility that
the degree $n$ of the polynomials $p_n$ only takes the values $n=0,1,\ldots,N$
and that the orthogonality is with respect to a positive measure having support
on a set of $N+1$ points.
Hahn \ref{Hah49} studied the $q$-analogue of this classification
(cf.\ \S2.5)
and he pointed out how the polynomials satisfying \eqtag{2.90} come out
as limit cases for $q\uparrow1$ of his classification.

The generic case for this classification is given by the {\sl Hahn polynomials}
$$
\eqalignno{
Q_n(x;\al,\be,N)
:=&
{}_3F_2\left[{-n,n+\al+\be+1,-x\atop\al+1,-N};1\right]
\cr
=&
\sum_{k=0}^n{(-n)_k\,(n+\al+\be+1)_k\,(-x)_k
\over
(\al+1)_k\,(-N)_k\,k!}\,.&\eq{2.100}
\cr}
$$
Here we assume $n=0,1,\ldots,N$ and we assume 
the notational convention that 
$$
{}_r F_s\left[{-n,a_2,\ldots,a_r\atop b_1,\ldots,b_s};z\right]
:=
\sum_{k=0}^n
{(-n)_k\,(a_2)_k\ldots(a_r)_k
\over
(b_1)_k\ldots (b_s)_k\,k!}\,z^k,\quad n\in\Zplus,
\eqno\eq{2.105}
$$
remains well-defined when some of the $b_i$ are possibly
non-positive integer but $\le-n$.
For the notation of Hahn polynomials and other families to be discussed
in this subsection we keep to the notation of
Askey \& Wilson \reff{AsWi85}{Appendix} and
Labelle's poster \ref{Lab90}.
There one can also find further references.

Hahn polynomials satisfy orthogonality relations
$$
\sum_{x=0}^N Q_n(x)\,Q_m(x)\,\rho(x)=\de_{n,m}\,{1\over\pi_n}\,,
\quad n,m=0,1\ldots,N,
\eqno\eq{2.110}
$$
where
$$
\eqalignno{
\rho(x)&:={N\choose x}{(\al+1)_x\,(\be+1)_{N-x}\over(\al+\be+2)_N}
&\eq{2.120}
\cr
\noalign{\hbox{and}}
\pi_n&:={N\choose n}{2n+\al+\be+1\over\al+\be+1}\,
{(\al+1)_n\,(\al+\be+1)_n\over(\be+1)_n\,(N+\al+\be+2)_n}\,.&
\eq{2.121}
\cr}
$$
We get positive weights $\rho(x)$ (or weights of fixed sign) if
$\al,\be\in(-1,\infty)\cup(-\infty,-N)$.

The other orthogonal polynomials coming out of this classification are
limit cases of Hahn polynomials. We get the following families.
\sLP
{\sl Krawtchouk polynomials}
$$
K_n(x;p,N):={}_2F_1(-n,-x;-N;p^{-1}),\quad n=0,1,\ldots,N,\quad 0<p<1,
\eqno\eq{2.125}
$$
with orthogonality measure having weights
$x\mapsto {N\choose x}\,p^x\,(1-p)^{N-x}$ on $\{0,1,\ldots,N\}$;
\sLP
{\sl Meixner polynomials}
$$
M_n(x;\be,c):={}_2F_1(-n,-x;\be;1-c^{-1}),\quad
0<c<1,\;\be>0,
\eqno\eq{2.126}
$$
with orthogonality measure having weights
$x\mapsto(\be)_x\,c^x/x!$ on $\Zplus$;
\sLP
{\sl Charlier polynomials}
$$
C_n(x;a):={}_2F_0(-n,-x;-;-a^{-1}),\quad a>0,
$$
with orthogonality measure having weights
$x\mapsto a^x/x!$ on $\Zplus$.

Note that Krawtchouk polynomials are Meixner polynomials \eqtag{2.126}
with $\be:=q^{-N}$.
Kraw\-tchouk and Meixner polynomials are limits of Hahn
polynomials, while Charlier polynomials are limits of both Krawtchouk
and Meixner polynomials.
In a certain sense, the very classical orthogonal polynomials are also
contained in the class discussed here, since Jacobi, Laguerre and Hermite
polynomials are limits of Hahn, Meixner and Charlier polynomials,
respectively. For instance,
$$
P_n^{(\al,\be)}(1-2x)=
{(\al+1)_n\over n!}\,\lim_{N\to\infty}
Q_n(Nx;\al,\be,N).
$$
See Table 1 (or rather a part of it) for a pictorial representation
of these families and their limit transitions.

The Krawtchouk, Meixner and Charlier polynomials are {\sl self-dual}, i.e.,
they satisfy
$$
p_n(x)=p_x(n),\quad x,n\in\Zplus\hbox{ or $x,n\in\{0,1,\ldots,N\}$.}
$$
Thus the orthogonality relations and dual orthogonality relations of these
polynomials essentially coincide
and the second order difference equation \eqtag{2.90} becomes the three
term recurrence relation after interchange of $x$ and $n$.
Then it can be arranged that the eigenvalue $\la_n$ in \eqtag{2.90}
becomes $n$.
By the self-duality, the $L^2$ completeness of the systems in case $N=\infty$
is also clear.

%\topinsert
%\vskip 13truecm
%\endinsert

For Hahn polynomials we have no self-duality, so the dual orthogonal system
will be different.
Observe that, in \eqtag{2.100}, we can write
$$
(-n)_k\,(n+\al+\be+1)_k=
\prod_{j=0}^{k-1}\bigl(-n(n+\al+\be+1)+j(j+\al+\be+1)\bigr),
$$
which is a polynomial of degree $k$ in $n(n+\al+\be+1)$.
Now define
$$
R_n(x(x+\al+\be+1);\al,\be,N):=
Q_x(n;\al,\be,N),\quad n,x=0,1,\ldots,N.
$$
Then $R_n$ ($n=0,1,\ldots,N$) extends to a polynomial of degree $n$:
$$
\eqalignno{
R_n(x(x+\al+\be+1);\al,\be,N)
&=
{}_3F_2\left[{-n,-x,x+\al+\be+1\atop\al+1,-N};1\right]
\cr
&=
\sum_{k=0}^n{(-n)_k\,(-x)_k\,(x+\al+\be+1)_k\over
(\al+1)_k\,(-N)_k\,k!}\,.
\cr}
$$
These are called the {\sl dual Hahn polynomials}.
The dual orthogonality relations implied by \eqtag{2.110} are the
orthogonality relations for the dual Hahn polynomials:
$$
\sum_{x=0}^N R_m(x(x+\al+\be+1))\,R_n(x(x+\al+\be+1))\,\pi_x
=
\de_{m,n}\,{1\over\rho(n)}\,.
$$
Here $\pi_x$ and $\rho(n)$ are as in \eqtag{2.120} and \eqtag{2.121}.
Thus the dual Hahn polynomials are also orthogonal polynomials.
The three term recurrence relation for Hahn polynomials translates as
a second order difference equation which is a slight generalization of
\eqtag{2.90}. It has the form
$$
a(x)\,p_n(\la(x+1))+b(x)\,p_n(\la(x))+c(x)\,p_n(\la(x-1))
=\la_n\,p_n(\la(x)),
\eqno\eq{2.130}
$$
where $\la(x):=x(x+\al+\be+1)$ is a quadratic function of $x$.
The Krawtchouk and Meixner polynomials are also limit cases of
the dual Hahn polynomials.

It is a natural question to ask for other orthogonal polynomials
being eigenfunctions of a second order difference equation of the form
\eqtag{2.130}.
A four-parameter family with this property are the
Racah polynomials, essentially known for a long time to the physicists
as Racah coefficients, which occur
in connection with
threefold tensor products of irreducible \rep s of the group $SU(2)$.
However, it was not recognized before the late seventies
(cf.\ Wilson \ref{Wil80})
that orthogonal polynomials are hidden in the Racah coefficients.
{\sl Racah polynomials} are defined by
$$
R_n(x(x+\ga+\de+1);\al,\be,\ga,\de)
:=
{}_4F_3\left[{-n,n+\al+\be+1,-x,x+\ga+\de+1
\atop
\al+1,\be+\de+1,\ga+1};1\right],
\eqno\eq{2.140}
$$
where
$\al+1$ or $\be+\de+1$ or $\ga+1=-N$ for some $N\in\Zplus$,
and where $n=0,1,\ldots,N$.
Similarly as for the dual Hahn polynomials it can be seen that
the Racah polynomial $R_n$ is indeed a polynomial of degree $n$ in
$\la(x):=x(x+\ga+\de+1)$.
The Racah polynomials are orthogonal \wrt weights $w(x)$ on the
points $\la(x)$, $x=0,1\ldots,N$, given by
$$
w(x):=
{(\ga+\de+1)_x\,((\ga+\de+3)/2)_x\,(\al+1)_x\,(\be+\de+1)_x\,(\ga+1)_x
\over
x!\,((\ga+\de+1)/2)_x\,(\ga+\de-\al+1)_x\,(\ga-\be+1)_x\,(\de+1)_x}\,.
$$
It is evident from \eqtag{2.140} that dual Racah polynomials are
again Racah polynomials with $\al,\be$ interchanged with $\ga,\de$.
Hahn and dual Hahn polynomials can be obtained as limit cases of
Racah polynomials.

Each orthogonality relation for the Racah polynomials is an explicit
evaluation of a finite sum. It is possible to interpret this finite
sum as a sum of residues coming from a contour integral in the complex plane.
This contour integral can also be considered and evaluated
for values of $N$ not necessarily in $\Zplus$.
For suitable values of the parameters the contour integral can then be deformed
to an integral over the imaginary axis and it gives rise to the
orthogonality relations for the {\sl Wilson polynomials\/}
(cf.\ Wilson \ref{Wil80}) defined by
$$
\displaylines{
\quad W_n(x^2;a,b,c,d)\hfill
\cr
\hfill:=(a+b)_n\,(a+c)_n\,(a+d)_n\,
{}_4F_3\left[{-n,n+a+b+c+d-1,a+ix,a-ix\atop
a+b,a+c,a+d};1\right].
\quad\eq{2.150}
\cr}
$$
Apparently, the \RHS\ defines a polynomial of degree $n$ in $x^2$.
If $a,b,c,d$ have positive real parts and complex parameters appear in
conjugate pairs then the functions $x\mapsto W_n(x^2)$ ($n\in\Zplus$)
of \eqtag{2.150} are orthogonal \wrt the measure $w(x)dx$ on $[0,\infty)$,
where
$$
w(x):=\left|
{\Ga(a+ix)\,\Ga(b+ix)\,\Ga(c+ix)\,\Ga(d+ix)\over\Ga(2ix)}\right|^2.
$$
The normalization in \eqtag{2.150} is such that the Wilson polynomials
are symmetric in their four parameters $a,b,c,d$.

The Wilson polynomials satisfy an eigenfunction equation of the form
$$
a(x)\,W_n((x+i)^2)+b(x)\,W_n(x^2)+c(x)\,W_n((x-i)^2)=\la_n\,W_n(x^2).
$$
So we have the new phenomenon that the difference operator at the \LHS\
shifts into the complex plane, out of the real interval on which the
Wilson polynomials are orthogonal.
This difference operator can be factorized as a product of two shift operators
of similar type. They have properties and applications analogous to the
shift operators for Jacobi polynomials discussed at the end of
\S\the\sectionnumber.1.
They also reduce the evaluation of the quadratic norms to
the evaluation of the integral of the weight function, but this last problem
is now much less trivial than in the Jacobi case.

Now, we can descend from the Wilson polynomials by limit transitions,
just as we did from the Racah polynomials.
On the ${}_3F_2$ level we thus get continuous analogues of the Hahn and dual
Hahn polynomials as follows.
\sLP
{\sl Continuous dual Hahn polynomials}:
$$
S_n(x^2;a,b,c):=(a+b)_n\,(a+c)_n\,{}_3F_2\left[
{-n,a+ix,a-ix\atop a+b,a+c};1\right],
$$
where $a,b,c$ have positive real parts; if one of these parameters is not real
then one of the other parameters is its complex conjugate.
The functions $x\mapsto S_n(x^2)$ are orthogonal with respect to the measure
$w(x)dx$ on $[0,\infty)$, where
$$
w(x):=\left|{\Ga(a+ix)\,\Ga(b+ix)\,\Ga(c+ix)\over\Ga(2ix)}\right|^2.
$$
\sLP
{\sl Continuous Hahn polynomials}:
$$
p_n(x;a,b,\bar a,\bar b):=i^n\,{(a+\bar a)_n\,(a+\bar b)_n\over n!}\,
{}_3F_2\left[{-n,n+a+\bar a+b+\bar b-1,a+ix\atop a+\bar a,a+\bar b};1\right],
$$
where $a,b$ have positive real part.
The polynomials $p_n$ are orthogonal on $\RR$ with respect to the measure
$|\Ga(a+ix)\,\Ga(b+ix)|^2\,dx$.
In Askey \& Wilson \reff{AsWi85}{Appendix}
only the symmetric case ($a,b>0$ or $a=\bar b$) of these polynomials occurs.
The general case was discovered by
Atakishiyev \& Suslov \ref{AtSu85}.

\centerline{\epsfxsize 375pt\epsfbox{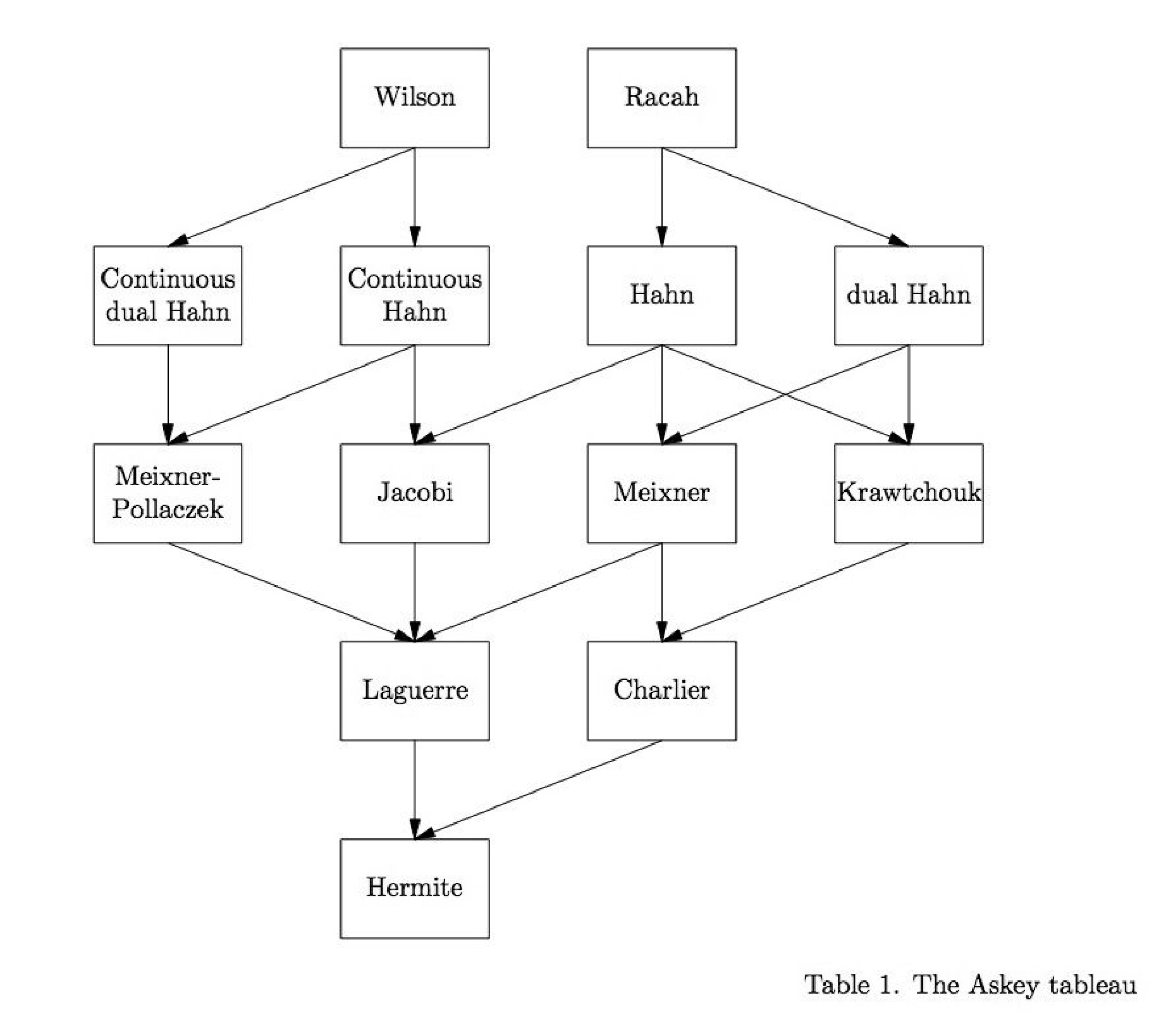}}

Jacobi polynomials are limit cases of continuous Hahn polynomials and
also directly of Wilson polynomials (with one pair of complex conjugate
parameters).
There is one further class of orthogonal polynomials on the ${}_2F_1$
level occurring as limit cases of orthogonal polynomials on the
${}_3F_2$ level:
\sLP
{\sl Meixner-Pollaczek polynomials}:
$$
P_n^{(a)}(x;\phi):={(2a)_n\over n!}\,
e^{in\phi}\,{}_2F_1(-n,a+ix;2a;1-e^{-2i\phi}),
\quad a>0,\;0<\phi<\pi.
$$
(Here we have chosen the normalization as in Labelle \ref{Lab90},
which is the same as in Pollaczek's original 1950 paper.)
They are orthogonal on $\RR$ with respect to the weight function
$x\mapsto e^{(2\phi-\pi)x}\,|\Ga(a+ix)|^2$. They can be considered
as continuous analogues of the Meixner and Krawtchouk polynomials.
They are limits of both continuous Hahn and continuous dual Hahn
polynomials. Laguerre polynomials are limit cases of Meixner-Pollaczek
polynomials.
Note that the last three families are analytic continuations,
both in $x$ and in the parameters, of
dual Hahn polynomials, Hahn polynomials and Meixner polynomials, respectively.

All families of orthogonal polynomials discussed until now,
together with the limit transitions between them, form
the {\sl Askey tableau\/} (or {\sl scheme} or {\sl chart})
{\sl of hypergeometric
orthogonal polynomials}. See Askey \& Wilson \reff{AsWi85}{Appendix},
Labelle \ref{Lab90} or Table 1.
See also \ref{Koo88} for group theoretic interpretations.

\Subsec3 {Big $q$-Jacobi polynomials}
These polynomials were hinted at by Hahn \ref{Hah49}
and explicitly introduced by
Andrews \& Askey \ref{AnAs85}.
Here we will show how their basic properties can be derived from a suitable
pair of shift operators.
We keep the convention of \S1 that $0<q<1$.

First we introduce $q$-integration by parts.
This will involve {\sl backward\/} and {\sl forward $q$-derivatives}:
$$
(D_q^-f)(x):={f(x)-f(qx)\over(1-q)x},\quad
(D_q^+f)(x):={f(q^{-1}x)-f(x)\over(1-q)x}\,.
$$
Here $D_q^-$ coincides with $D_q$ introduced in \eqtag{1.36}.

\Prop{2.10}
If $f$ and $g$ are continuous on $[-d,c]$ ($c,d\ge0$) then
$$
\int_{-d}^c(D_q^-f)(x)\,g(x)\,d_qx
=
f(c)\,g(q^{-1}c)-f(-d)\,g(-q^{-1}d)-
\int_{-d}^c f(x)\,(D_q^+g)(x)\,d_qx.
$$

\Proof
$$
\eqalignno{
\int_0^c
(D_q^-f)&(x)\,g(x)\,d_q(x)
=\sum_{k=0}^\infty\bigl(f(cq^k)-f(cq^{k+1})\bigr)\,g(cq^k)
\cr
&=
\lim_{N\to\infty}\biggl\{
f(c)\,g(q^{-1}c)-f(cq^{N+1})\,g(cq^N)+
\sum_{k=0}^N f(cq^k)\,\bigl(g(cq^k)-g(cq^{k-1})\bigr)\biggr\}
\cr
&=
f(c)\,g(q^{-1}c)-f(0)\,g(0)+
\sum_{k=0}^\infty f(cq^k)\,\bigl(g(cq^k)-g(cq^{k-1})\bigr)
\cr
&=
f(c)\,g(q^{-1}c)-f(0)\,g(0)-
\int_0^c f(x)\,(D_q^+g)(x)\,d_qx.
\cr}
$$
Now apply \eqtag{1.50}.\hhalmos

\bPP
Let
$$
w(x;a,b,c,d;q):=
{(qx/c,-qx/d;q)_\infty\over
(qax/c,-qbx/d;q)_\infty}\,.
\eqno\eq{2.198}
$$
Note that $w(x;a,b,c,d;q)>0$ on $[-d,c]$ if $c,d>0$ and
$$
\left[\,-{c\over dq}<a<{1\over q}\;\&\;
-{d\over cq}<b<{1\over q}\,\right]\quad {\rm or}\quad
\bigl[a=c\al\;\&\;b=-d\bar\al\hbox{ for some $\al\in\CC\backslash\RR$.}\bigr]
\eqno\eq{2.199}
$$
{}From now on we assume that these inequalities hold.
If $f$ is continuous on $[-d,c]$ then
$$
\lim_{q\uparrow1}\int_{-d}^c
f(x)\,w(x;q^\al,q^\be,c,d;q)\,d_qx=
\int_{-d}^c f(x)\,(1-c^{-1}x)^\al\,(1+d^{-1}x)^\be\,dx,
\eqno\eq{2.200}
$$
Thus the measure $w(x;a,b,c,d;q)\,d_qx$ can be considered as a
$q$-analogue of the Jacobi polynomial orthogonality measure shifted
to an arbitrary finite interval containing 0.

Since
$$
w(q^{-1}c;a,b,c,d;q)=0=w(-q^{-1}d;a,b,c,d;q),
$$
we get from Proposition \thtag{2.10} that for any two polynomials
$f,g$ the following holds.
$$
\displaylines{
\quad \int_{-d}^c (D_q^-f)(x)\,g(x)\,w(x;qa,qb,c,d;q)\,d_qx\hfill
\cr
\hfill=-\int_{-d}^c f(x)\,
\bigl[D_q^+\bigl(g(.)\,w(.;qa,qb,c,d;q)\bigr)\bigr](x)\,d_qx.\quad
\eq{2.210}
\cr}
$$

Define
$$
\eqalignno{
(D_q^{+,a,b}f)(x):=&
{\bigl[D_q^+\bigl(w(.;qa,qb,c,d;q)\,f(.)\bigr)\bigr](x)\over w(x;a,b,c,d;q)}
\cr
=&
(1-q)^{-1}\,x^{-1}(1-x/c)\,(1+x/d)\,f(q^{-1}x)
\cr
&\quad
-(1-q)^{-1}\,x^{-1}\,(1-qax/c)\,(1+qbx/d)\,f(x).
&\eq{2.215}
\cr}
$$
In the following $D_q^-f(x)$ or $D_q^-(f(x))$ will mean $(D_q^-f)(x)$.
Also $D_q^{+,a,b}f(x)$ or $D_q^{+,a,b}(f(x))$ will mean
$(D_q^{+,a,b}f)(x)$.
A simple computation yields the following.
$$
\displaylines{
\hfill D_q^-\,x^n={1-q^n\over1-q}\,x^{n-1},\hfill\eq{2.220}
\cr
\qquad D_q^{+,a,b}\bigl((q^2ax/c;q)_{n-1}\bigr)\hfill
\cr
\qquad={q^{-n}a^{-1}-qb\over(1-q)d}\,(qax/c;q)_n+
\hbox{terms of degree $\le n-1$ in $x$}.\hfill\eq{2.230}
\cr}
$$

Define the {\sl big $q$-Jacobi polynomial\/}
$\Pt_n(x;a,b,c,d;q)$
as the monic orthogonal polynomial of degree $n$ in $x$ \wrt
the measure $w(x;a,b,c,d;q)\,d_qx$ on $[-d,c]$.
Later we will introduce another normalization for these polynomials and
we will
then write them as $P_n$ instead of $\Pt_n$.
It follows from \eqtag{2.200} that
$$
\lim_{q\uparrow1}\Pt_n(x;q^\al,q^\be,c,d;q)=
\const\,P_n^{(\al,\be)}\left({d-c+2x\over d+c}\right),
$$
a Jacobi polynomial shifted to the interval $[-d,c]$.
Big $q$-Jacobi polynomials with $d=0$, $c=1$ are called
{\sl little $q$-Jacobi polynomials}.

Two simple consequences of the definition are:
$$
\eqalignno{
\Pt_n(-x;a,b,c,d;q)&=(-1)^n\,\Pt_n(x;b,a,d,c;q)
&\eq{2.235}
\cr
\noalign{\hbox{and}}
\Pt_n(\la x;a,b,c,d;q)&=\la^n\,\Pt_n(x;a,b,\la^{-1}c,\la^{-1}d;q),\quad
\la>0.
\cr}
$$
It follows from \eqtag{2.210}, \eqtag{2.220} and \eqtag{2.230} that
$D_q^-$ and $D_q^{+,a,b}$ act as shift operators on the big $q$-Jacobi
polynomials:
$$
\eqalignno{
D_q^-\,\Pt_n(x;a,b,c,d;q)&=
{1-q^n\over1-q}\,\Pt_{n-1}(x;qa,qb,c,d;q),
&\eq{2.240}
\cr
D_q^{+,a,b}\,\Pt_{n-1}(x;qa,qb,c,d;q)&=
{q^2ab-q^{-n+1}\over(1-q)cd}\,\Pt_n(x;a,b,c,d;q).
&\eq{2.250}
\cr}
$$
Composition of \eqtag{2.240} and \eqtag{2.250} yields a second order
$q$-difference equation for the big $q$-Jacobi polynomials:
$$
D_q^{+,a,b}\,D_q^-\,\Pt_n(x;a,b,c,d;q)=
{q(1-q^{-n})(1-q^{n+1}ab)\over(1-q)^2cd}\,\Pt_n(x;a,b,c,d;q).
\eqno\eq{2.260}
$$
The \LHS\ of \eqtag{2.260} can be rewritten as
$$
\displaylines{
\qquad
A(x)\,\bigl((D_q^-)^2\Pt_n\bigr)(q^{-1}x)+B(x)\,(D_q^-\Pt_n)(q^{-1}x)\hfill
\cr
\hfill=a(x)\Pt_n(q^{-1}x)+b(x)\Pt_n(x)+c(x)\Pt_n(qx)\qquad
\eq{2.270}
\cr}
$$
for certain polynomials $A,B$ and $a,b,c$.
Here
$$
\eqalignno{
A(x)&=
q^{-1}\,(1-qax/c)\,(1+qbx/d),
\cr
B(x)&=
{(1-qb)c-(1-qa)d-(1-q^2ab)x
\over
(1-q)cd}\,.
&\eq{2.272}
\cr}
$$
Compare \eqtag{2.260}, \eqtag{2.270} with \eqtag{2.5}, \eqtag{2.90}.
Apparently we have here another extension of the concept of classical
orthogonal polynomials.
Two other properties point into the same direction.
First, by \eqtag{2.240} the polynomials
$D_q^-\Pt_{n+1}$ ($n=0,1,2,\ldots$) form again a system of orthogonal
polynomials.
Second, iteration of \eqtag{2.250} yields a Rodrigues type formula.

It follows from \eqtag{2.215} that
$$
(D_q^{+,a,b}f)\left({c\over qa}\right)=
-{(1-qa)\,(1+qad/c)\over qad(1-q)}\,
f\left({c\over q^2a}\right).
$$
Combination with \eqtag{2.250} yields the recurrence
$$
\displaylines{
\qquad{q^2ab-q^{-n+1}\over(1-q)cd}\,\Pt_n\left({c\over qa};a,b,c,d;q\right)
\hfill\cr
\hfill=-{(1-qa)(1+qad/c)\over qad(1-q)}\,
\Pt_{n-1}\left({c\over q^2a};qa,qb,c,d;q\right).\qquad
\cr}
$$
By iteration we get an evaluation of the big $q$-Jacobi polynomial at
a special point:
$$
\Pt_n\left({c\over qa};a,b,c,d;q\right)
=\left(c\over qa\right)^n\,
{(qa;q)_n\,(-qad/c;q)_n\over(q^{n+1}ab;q)_n}\,.
\eqno\eq{2.278}
$$
Now we will normalize the big $q$-Jacobi polynomials such that
they take the value 1 at $c/(qa)$:
$$
P_n(x;a,b,c,d;q):={\Pt_n(x;a,b,c,d;q)\over
\Pt_n(c/(qa);a,b,c,d;q)}\,.
\eqno\eq{2.280}
$$
Note that the value at $-d/(qb)$ now follows from \eqtag{2.235} and
\eqtag{2.278}:
$$
P_n\left({-d\over qb};a,b,c,d;q\right)=
\left(- {ad\over bc}\right)^n\,
{(qb;q)_n\,(-qbc/d;q)_n\over(qa;q)_n\,(-qad/c;q)_n}\,.
\eqno\eq{2.281}
$$
It follows from \eqtag{2.235} and \eqtag{2.281}
that
$$P_n(-x;a,b,c,d;q)=
\left(- {ad\over bc}\right)^n\,
{(qb;q)_n\,(-qbc/d;q)_n\over(qa;q)_n\,(-qad/c;q)_n}\,
P_n(x;b,a,d,c;q).
$$
The following lemma, proved by use of \eqtag{1.37}, gives the $q$-analogue
of a Taylor series expansion.

\Lemma{2.20}
If $f(x):=\sum_{k=0}^n c_k\,(qax/c;q)_k$ then
$$
\bigl((D_q^-)^kf\bigr)\left({c\over q^{k+1}a}\right)=
c_k\,(-1)^k\,\left({qa\over c}\right)^k\,q^{k(k-1)/2}\,
{(q;q)_k\over(1-q)^k}\,.
$$

\bPP
Now put $f(x):=P_n(x;a,b,c,d;q)$ in Lemma \thtag{2.20} and substitute
\eqtag{2.240} (iterated) and \eqtag{2.280}.
Then the $c_k$ can be found explicitly and we obtain a representation by
a $q$-hypergeometric series:
$$
\eqalignno{
P_n(x;a,b,c,d;q)&=
\sum_{k=0}^n {(q^{-n};q)_k\,(q^{n+1}ab;q)_k\,(qax/c;q)_k\,q^k\over
(qa;q)_k\,(-qad/c;q)_k\,(q;q)_k}
\cr
&=
{}_3\phi_2\left[{q^{-n},q^{n+1}ab,qax/c\atop qa,-qad/c};q,q\right].
&\eq{2.285}
\cr}
$$
Andrews \& Askey \reff{AnAs85}{(3,28)}
use the notation $P_n^{(\al,\be)}(x;c,d:q)$,
which coincides, up to a constant factor, with our
$p_n(x;q^\al,q^\be,c,d;q)$.

It follows by combination of \eqtag{2.235} and \eqtag{2.285} that
$$
{P_n(x;a,b,c,d;q)\over
P_n(-d/(qb);a,b,c,d;q)}=
{}_3\phi_2\left[{q^{-n},q^{n+1}ab,-qbx/d\atop qb,-qbc/d};q,q\right],
\eqno\eq{2.286}
$$
where the denominator of the \LHS\ is explicitly given by \eqtag{2.281}.

Next we will compute the quadratic norms.
By \eqtag{2.210}, \eqtag{2.240} and \eqtag{2.250} we get the recurrence
$$
\displaylines{
\hfill\int_{-d}^c(\Pt_n^2w)(x;a,b,c,d;q)\,d_qx
=
{q^{n-1}\,(1-q^n)\,cd\over1-q^{n+1}ab}\,
\int_{-d}^c(\Pt_{n-1}^2w)(x;qa,qb,c,d;q)\,d_qx\hfill
\cr
\hfill={q^{n(n-1)/2}\,(q;q)_n\,(cd)^n\over(q^{n+1}ab;q)_n}\,
\int_{-d}^c w(x;q^na,q^nb,c,d;q)\,d_qx,
\quad\eq{2.290}
\cr}
$$
where the second equality follows by iteration.
Now the $q$-integral of $w$ from $-d$ to $c$ can be rewritten as a sum
of two ${}_2\phi_1$'s of argument $q$ of the form of the \LHS\ of
\eqtag{1.220}. Evaluation by \eqtag{1.220} yields
$$
\int_{-d}^c w(x;a,b,c,d;q)\,d_qx=
(1-q)c\,{(q,-d/c,-qc/d,q^2ab;q)_\infty\over
(qa,qb,-qbc/d,-qad/c;q)_\infty}\,.
\eqno\eq{2.294}
$$
Together with \eqtag{2.290} this yields:
$$
{\int_{-d}^c(\Pt_n^2w)(x;a,b,c,d;q)\,d_qx\over
\int_{-d}^c w(x;a,b,c,d;q)\,d_qx}
=q^{n(n-1)/2}\,(cd)^n\,
{(q,qa,qb,-qbc/d,-qad/c;q)_n\over
(q^2ab;q)_{2n}\,(q^{n+1}ab;q)_n}\,.
\eqno\eq{2.295}
$$

Now we can compute the coefficients in the three term recurrence relation
for the big $q$-Jacobi polynomials $\Pt_n(x):=\Pt_n(x;a,b,c,d;q)$:
$$
x\,\Pt_n(x)=\Pt_{n+1}(x)+B_n\,\Pt_n(x)+C_n\,\Pt_{n-1}(x).
\eqno\eq{2.296}
$$
Then $C_n$ is the quotient of the \RHS s of \eqtag{2.295} for degree $n$ and
for degree $n-1$, respectively, so
$$
C_n=
{q^{n-1}\,(1-q^n)\,(1-q^na)\,(1-q^n b)\,(1-q^n ab)\,(d+q^n bc)\,(c+q^n ad)
\over
(1-q^{2n-1}ab)\,(1-q^{2n}ab)^2\,(1-q^{2n+1}ab)}\,.
\eqno\eq{2.297}
$$
Then we obtain $B_n$ by substitution of $x:=c/(qa)$ in \eqtag{2.296}, in view of
\eqtag{2.278}.

From \eqtag{2.285} one sees that $P_n(x;a,b,c,d;q)$ is homogeneous
of degree 0 in the three variables $x,c,d$. Therefore, only three
of the four parameters $a,b,c,d$ are essential. Probably for this reason,
Gasper \& Rahman \reff{GaRa90}{(7.3.10)}
use the following notation for big $q$-Jacobi polynomials:
$$
P_n^{\rm GR}(x;a,b,c;q)=P_n(x;a,b,c;q):=
{}_3\phi_2\left[{q^{-n},abq^{n+1},x\atop aq,cq};q,q\right].
$$
Their notation is related to the notation \eqtag{2.285} in the present paper
by
$$
\eqalignno{
P_n(x;a,b,c,d;q)&=P_n^{\rm GR}(ac^{-1}qx;a,b,-ac^{-1}d;q),
\cr
P_n^{\rm GR}(x;q,b,c;q)&=P_n(x;a,b,aq,-cq;q).
\cr}
$$

\Subsec4 {Little $q$-Jacobi polynomials}
These polynomials are the most straightforward $q$-analogues of the
Jacobi polynomials.
They were first observed by Hahn \ref{Hah49}
and studied in more detail by Andrews \& Askey \ref{AnAs77}.
Here we will study them as a special case of the big $q$-Jacobi polynomials.

When we specialize the big $q$-Jacobi polynomials
\eqtag{2.280} to $c=1$, $d=0$ and normalize them such that they take the value
1 at 0 then we obtain the {\sl little $q$-Jacobi polynomials\/}
$$
\eqalignno{
p_n(x;a,b;q)
:=&
{P_n(x;b,a,1,0;q)\over P_n(0;b,a,1,0;q)}
\cr
=&
{}_2\phi_1(q^{-n},q^{n+1}ab;qa;q,qx)
&\eq{2.300}
\cr
=&
(-qb)^{-n}\,q^{-n(n-1)/2}\,{(qb;q)_n\over(qa;q)_n}\,
{}_3\phi_2\left[{q^{-n},q^{n+1}ab,qbx\atop qb,0};q,q\right].
&\eq{2.310}
\cr}
$$
Here \eqtag{2.300} follows by letting $d\to0$ in \eqtag{2.286},
while \eqtag{2.310} follows from \eqtag{2.285} and \eqtag{2.281}.
With the equality of \eqtag{2.300} and \eqtag{2.310} we have reobtained
the transformation formula \eqtag{1.204} for terminating ${}_2\phi_1$ series.
{}From \eqtag{2.310} we obtain
$$
p_n(q^{-1}b^{-1};a,b;q)=
(-qb)^{-n}\,q^{-n(n-1)/2}\,{(qb;q)_n\over(qa;q)_n}\,.
$$
{}From \eqtag{2.300} and \eqtag{1.202} we obtain
$$
p_n(1;a,b;q):=(-a)^n\,q^{n(n+1)/2}\,{(qb;q)_n\over(qa;q)_n}\,.
$$

Little $q$-Jacobi polynomials satisfy the orthogonality relations
$$
\displaylines{
\quad{1\over B_q(\al+1,\be+1)}\,
\int_0^1 p_n(t;q^\al,q^\be;q)\,
p_m(t;q^\al,q^\be;q)\,t^\al\,
{(qt;q)_\infty\over(q^{\be+1}t;q)_\infty}\,d_qt\hfill
\cr
\hfill
=\de_{n,m}\,
{q^{n(\al+1)}\,(1-q^{\al+\be+1})\,(q^{\be+1};q)_n\,(q;q)_n
\over
(1-q^{2n+\al+\be+1})\,(q^{\al+1};q)_n\,(q^{\al+\be+1};q)_n}\,.\quad
\eq{2.320}
\cr}
$$
The orthogonality measure is the measure of the $q$-beta integral
\eqtag{1.153}. 
For positivity and convergence we require that $0<a<1/q$, $b<1/q$ (after we have
replaced $q^\al$ by $a$ and $q^\be$ by $b$ in \eqtag{2.320}).
It is maybe not immediately seen that \eqtag{2.320} is a limit case of
\eqtag{2.295}. However, we can establish this by observing the weak convergence
as $d\downarrow 0$ of the  normalized measure
const.$\times w(x;q^\be,q^\al,1,d;q)d_qx$ on $[-d,1]$ to the normalized measure
in \eqtag{2.320} on [0,1]:
$$
\displaylines{
\quad\lim_{d\downarrow0}
{\int_{-d}^1 (q^{\be+1}t;q)_n\,w(t;q^\be,q^\al,1,d;q)\,d_qt
\over
\int_{-d}^1w(t;q^\be,q^\al,1,d;q)\,d_qt}
=\lim_{d\downarrow0}{(q^{\be+1};q)_n\,(-q^{\be+1}d;q)_n\over
(q^{\al+\be+2};q)_n}\hfill
\cr
\hfill={(q^{\be+1};q)_n\over(q^{\al+\be+2};q)_n}
={1\over B_q(\al+1,\be+1)}\,\int_0^1(q^{\be+1}t;q)_n\,
{(qt;q)_\infty\over(q^{\be+1}t;q)_\infty}\,t^\al\,d_qt.\quad
\cr}
$$
Here the first equality follows from \eqtag{2.294} and \eqtag{2.198},
while the last equality follows from \eqtag{1.153}.

Little $q$-Jacobi polynomials are the $q$-analogues of the Jacobi polynomials
shifted to the interval $[0,1]$:
$$
\lim_{q\uparrow1}p_n(x;q^\al,q^\be;q)=
{P_n^{(\al,\be)}(1-2x)\over P_n^{(\al,\be)}(1)}\,.
$$
If we fix $b<1$ (possibly 0 or negative) and put $a=q^\al$ then little
$q$-Jacobi polynomials tend for $q\uparrow1$ to Laguerre polynomials:
$$
\lim_{q\uparrow1}
p_n\left({(1-q)x\over1-b};q^\al,b;q\right)=
{L_n^\al(x)\over L_n^\al(0)}\,.
$$
In particular, little $q$-Jacobi polynomials
$p_n(x;a,0;q)$, called {\sl Wall polynomials}, are $q$-analogues of Laguerre
polynomials.

Analogous to the quadratic tansformations for
Jacobi polynomials (cf.\ \reff{Erd2}{10.9(4), (21) and (22)}) we can find
{\sl quadratic transformations\/} between little and big $q$-Jacobi
polynomials:
$$
\eqalignno{
P_{2n}(x;a,a,1,1;q)&=
{p_n(x^2;q^{-1},a^2;q^2)\over
p_n((qa)^{-2};q^{-1},a^2;q^2)}\,,
&\eq{2.330}
\cr
P_{2n+1}(x;a,a,1,1;q)&=
{x\,p_n(x^2;q,a^2;q^2)\over
(qa)^{-1}\,p_n((qa)^{-2};q,a^2;q^2)}\,.
&\eq{2.340}
\cr}
$$
The proof is also analogous: by use of the orthogonality properties and
normalization of the polynomials.

Just as Jacobi polynomials tend to Bessel functions by
$$
\displaylines{
\quad
{P_{n_N}^{(\al,\be)}(1-x^2/(2N^2))
\over
P_{n_N}^{(\al,\be)}(1)}
=
{}_2F_1\left(
-n_N,n_N+\al+\be+1;\al+1;
{x^2\over4N^2}
\right)
\hfill\cr\hfill
{\buildrel N\to\infty\over\longrightarrow}
{}_0F_1(-;\al+1;-(\la x/2)^2)
=
(\la x/2)^{-\al}\,\Ga(\al+1)\,J_\al(\la x),\quad
n_N/N\to\la\hbox{ as $N\to\infty$,}
\cr}
$$
little $q$-Jacobi polynomials tend to Jackson's third $q$-Bessel function
\eqtag{1.380}:
$$
\eqalignno{
\lim_{N\to\infty}
p_{N-n}(q^Nx;a,b;q)
&=
\lim_{N\to\infty}
{}_2\phi_1(q^{-N+n},abq^{N-n+1};aq;q,q^{N+1}x)
\cr
&=
{}_1\phi_1(0;aq;q,q^{n+1}x),
\cr}
$$
cf.\ Koornwinder \& Swarttouw \reff{KoSw}{Prop.\ A.1}.

\Subsec5 {Hahn's classification}
Hahn \ref{Hah49}
classified all families of orthogonal polynomials $p_n$
($n=0,1,\ldots,N$ or $n=0,1,\ldots$)
which are eigenfunctions of a second order $q$-difference equation, i.e.,
$$
\eqalignno{
A(x)\,(D_q^2\,p_n)(q^{-1}x)+B(x)\,(D_q\,p_n)(q^{-1}x)
&=
\la_n\,p_n(x)
&\eq{2.400}
\cr
\noalign{\hbox{or equivalently}}
a(x)\,p_n(q^{-1}x)+b(x)\,p_n(x)+c(x)\,p_n(qx)
&=
\la_n\,p_n(x),
\cr}
$$
where $A,B$ and $a,b,c$ are fixed polynomials.
Necessarily, $A$ is of degree $\le2$ and $B$ is of degree $\le1$.
The eigenvalues $\la_n$ will be completely determined by $A$ and $B$.
One distinguishes cases depending on the degrees of $A$ and $B$ and the
situation of the zeros of $A$.
For each case one finds a family of $q$-hypergeometric polynomials
satisfying \eqtag{2.400} which, moreover, satisfies an
explicit three term recurrence relation
$$
x\,p_n(x)=p_{n+1}(x)+B_n\,p_n(x)+C_n\,p_{n-1}(x).
$$
(Here we assumed the $p_n$ to be monic.)
Then we will have orthogonal polynomials \wrt a positive orthogonality measure
iff $C_n>0$ for all $n$.
A next problem is to find the explicit orthogonality measure. For a given
family of $q$-hypergeometric polynomials depending on parameters the type of
this measure may vary with the parameters.
Finally the limit transitions between the various families of orthogonal
polynomials can be examined.

In essence this program has been worked out by Hahn \ref{Hah49},
but his paper is somewhat sketchy in details.
Unfortunately, there is no later publication, where the details have all been
filled in.
In Table 2 we give a {\sl $q$-Hahn tableau}:
a $q$-analogue of that part of the Askey tableau (Table 1)
which is dominated by the Hahn polynomials.
In the ${}_r\phi_s$ formulas in the Table we have omitted the last but one
parameter denoting the base except when this is different from $q$.
The arrows denote limit transitions.
The box for case 3a, which occurred in an earlier version, is now suppressed,
because it does not correspond to a class of orthogonal polynomials.
Below we will give a brief discussion of each case.
We do not claim completeness.

\centerline{\epsfxsize 375pt\epsfbox{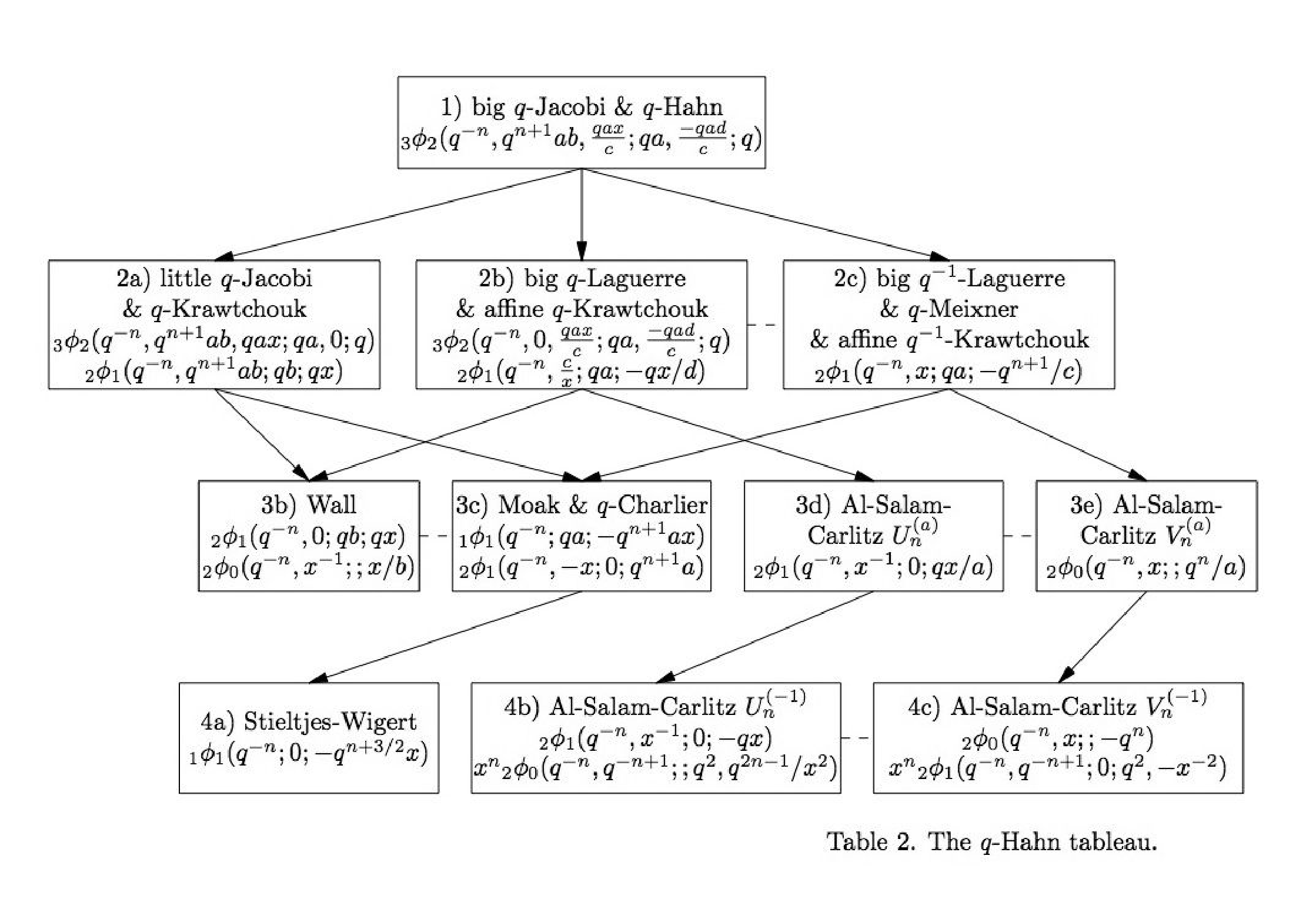}}

Several new phenomena occur here:
\LP
1. Within one class of $q$-hypergeometric polynomials we may obtain, depending
on the values of the parameters, either a $q$-analogue of the
Jacobi-Laguerre-Hermite class or of the (discrete)
Hahn-Krawtchouk-Meixner-Charlier class.
\LP
2. $q$-Analogues of Jacobi, Laguerre and Hermite polynomials occur in a
``little'' version (corresponding to Jacobi polynomials on $[0,c]$,
Laguerre polynomials on $[0,\infty)$ and Hermite polynomials which are even or
odd functions) and a ``big'' version ($q$-analogues of Jacobi, Laguerre
and Hermite polynomials of arbitrarily shifted argument).
\LP
3. There is a duality under the transformation $q\mapsto q^{-1}$,
denoted in Table 2 by dashed lines.
We insist on the convention $0<q<1$, but we can rewrite $q^{-1}$-hypergeometric
orthogonal polynomials as $q$-hypergeometric orthogonal polynomials
and thus sometimes obtain another family.
\LP
4. Cases may occur where the orthogonality measure is not uniquely determined.

%\topinsert
%\vskip 8.9truein
%\endinsert

\sPP
We now briefly discuss each case occuring in Table 2.

\sLP
1) {\sl Big $q$-Jacobi polynomials} $P_n(x;a,b,c,d;q)$,
defined by \eqtag{2.280}, form the generic case in this classification.
$A(x)$ and $B(x)$ in \eqtag{2.400} are given by \eqtag{2.272} and
$\la_n$ is as in the right hand side of \eqtag{2.260}.
{}From the explicit expression for $C_n$ in \eqtag{2.297} we get the values of
$a,b,c,d$ for which there is a positive orthogonality measure. For $c,d>0$
these are given by \eqtag{2.199}.

The {\sl $q$-Hahn polynomials\/} can be obtained as special big $q$-Jacobi
polynomials with $-qad/c=q^{-N}$, $n=0,1,\ldots,N$ ($N\in\Zplus$). For
convenience we may take $c=qa$. The $q$-Hahn polynomials are usually
(cf.\ \reff{GaRa90}{(7.2.21)})
notated as
$$
Q_n(x;a,b.N;q):=
{}_3\phi_2\left[
{q^{-n},abq^{n+1},x\atop aq,q^{-N}};q,q\right]
=
\sum_{k=0}^n
{(q^{-n};q)_k\,(abq^{n+1};q)_k\,(x;q)_k
\over
(aq;q)_k\,(q^{-N};q)_k\,(q;q)_k}\,q^k.
$$
Here the convention regarding lower parameters $q^{-N}$
($N\in\Zplus$) is similar to the convention for \eqtag{2.105}.

$q$-Hahn polynomials satisfy orthogonality relations
$$
\sum_{x=0}^N(Q_nQ_m)(q^{-x};a,b,N;q)\,
{(aq;q)_x\,(bq;q)_{N-x}
\over
(q;q)_x\,(q;q)_{N-x}}\,
(aq)^{-x}=0,\quad n\ne m.
$$
The weights are positive in one of the three following cases:
(i) $b<q^{-1}$ and $0<a<q^{-1}$;
(ii) $a,b>q^{-N}$;
(iii) $a<0$ and $b>q^{-N}$.
For $q\uparrow1$ the polynomials tend to ordinary Hahn polynomials \eqtag{2.100}:
$$
\lim_{q\uparrow1}Q_n(q^{-x};q^\al,q^\be,N;q)=
Q_n(x;\al,\be,N).
$$
\sLP
2a) {\sl Little $q$-Jacobi polynomials\/} $p_n(x;a,b;q)$,
given by \eqtag{2.300}, \eqtag{2.310},
were discussed in \S2.4.
When we put $b:=q^{-N-1}$ in \eqtag{2.310},
we obtain $q$-analogues of Krawtchouk polynomials \eqtag{2.125}:
the {\sl $q$-Krawtchouk polynomials}
$$
K_n(x;b,N;q)
:=
{}_3\phi_2\left[{q^{-n},-b^{-1}q^n,x\atop 0,q^{-N}};q,q\right]
=
\lim_{a\to 0}Q_n(x;a,-(qba)^{-1},N;q)
$$
with orthogonality relations
$$
\sum_{x=0}^N(K_nK_m)(q^{-x};b,N;q)\,
{(q^{-N};q)_x\over(q;q)_x}\,(-b)^x
=0,\quad n\ne m,\quad n,m=0,1,\ldots,N,
$$
and limit transition
$$
\lim_{q\uparrow1}K_n(q^{-x};b,N;q)=K_n(x;b/(b+1),N).
$$
See \reff{GaRa90}{Exercise 7.8(i)} and the reference given there.
Note that the $q$-Krawtchouk polynomials $K_n(q^{-x};b,N;q)$
are not self-dual under the
interchange of $x$ and $n$.
\sLP
2b) {\sl Big $q$-Laguerre polynomials}, these are big $q$-Jacobi
polynomials with $b=0$:
$$
\eqalignno{
P_n(x;a,0,c,d;q)
&=
{}_3\phi_2\left[{q^{-n},0,qax/c\atop qa,-qad/c};q,q\right]
&\eq{2.410}
\cr
&=
{1\over(-q^{-n}ca^{-1}d^{-1};q)_n}\,
{}_2\phi_1\left[{q^{-n},c/x\atop qa};q,-qx/d\right].
\cr}
$$
Here the second equality follows by \eqtag{1.206}.
Note that
$$
\lim_{q\uparrow1}P_n(x;q^\al,0,c,(1-q)^{-1};q)=
{}_1F_1(-n;\al+1;c-x)=\const L_n^\al(c-x).
$$

Another family of $q$-analogues of Krawtchouk polynomials, called
{\sl affine $q$-Kraw\-tchouk polynomials}, can be obtained from \eqtag{2.410}
by putting $-qad/c=q^{-N}$:
$$
K_n^{Aff}(x;a,N;q):=
{}_3\phi_2\left[{q^{-n},0,x\atop aq,q^{-N}};q,q\right]
=
Q_n(x;a,0,N;q)
$$
with orthogonality relations
$$
\sum_{x=0}^N(K_n^{Aff}K_m^{Aff})(q^{-x};a,N;q)\,
{(aq;q)_x\,(aq)^{-x}\over(q;q)_x\,(q;q)_{N-x}}=0,\quad
n\ne m,\quad n,m=0,1,\ldots,N,
$$
and limit transition
$$
\lim_{q\uparrow1} K_n^{Aff}(q^{-x};a,N;q)
=K_n(x;1-a,N).
$$
See \reff{GaRa90}{Exercise 7.11} and the references given there.

\sLP
2c) {\sl $q$-Meixner polynomials}
$$
M_n(x;a,c;q):={}_2\phi_1\left[{q^{-n},x\atop qa};q,{-q^{n+1}\over c}\right].
\eqno\eq{2.420}
$$
For certain values of the parameters
these can be considered as $q$-analogues of the Meixner polynomials
\eqtag{2.126}, cf.\ \reff{GaRa90}{Exercise 7.12}.
However, when we write these polynomials as
$$
M_n(qax;a,-\bar a^{-1};q)=
{}_2\phi_1\left[{q^{-n},qax\atop qa};q,q^{n+1}\bar a\right]
=\lim_{d\to\infty}P_n(x;a,-\bar a d,1,d;q)
$$
and $a$ is complex but not real, then these polynomials also become orthogonal
(not documented in the literature) and
they can be considered as $q$-analogues
of Laguerre polynomials of shifted argument.
Moreover, these polynomials can be obtained from the big $q$-Laguerre
polynomials \eqtag{2.410} by the transformation $q\mapsto q^{-1}$.
Therefore, we call these polynomials also
{\sl big $q^{-1}$-Laguerre polynomials}.

For $a:=q^{-N-1}$ the polynomials \eqtag{2.420} become yet another family
of $q$-analogues of the Krawtchouk polynomials, which we will call
{\sl affine $q^{-1}$-Krawtchouk polynomials},
because they are obtained from affine $q$-Krawtchouk polynomials by changing
$q$ into $q^{-1}$.
These polynomials, written as
$$
M_n(x;q^{-N-1},-b^{-1};q):=
{}_2\phi_1(q^{-n},x;q^{-N};q,bq^{n+1})
=
\lim_{a\to\infty}Q_n(x;a,b,N;q),
$$
where $n=0,1,\ldots,N$, have orthogonality relations
$$
\sum_{x=0}^N(M_nM_m)(q^{-x};q^{-N-1},-b^{-1};q)\,
{(bq;q)_{N-x}\,(-1)^{N-x}\,q^{x(x-1)/2}
\over
(q;q)_x\,(q;q)_{N-x}}=0,\quad n\ne m,
$$
and limit transition
$$
\lim_{q\uparrow1}M_n(q^{-x};q^{-N-1},-b^{-1};q)=
K_n(x;b^{-1},N).
$$
See Koornwinder \ref{Koo89}.

\sLP
3b) {\sl Wall polynomials\/}
(cf.\ Chihara \reff{Chi78}{\S VI.11})
are special little $q$-Jacobi polynomials
$$
\eqalignno{
p_n(x;a,0;q)&={}_2\phi_1(q^{-n},0;qa;q,qx)
\cr
&={1\over(q^na^{-1};q)_n}\,{}_2\phi_0(q^{-n},x^{-1};-;q,x/a).
\cr}
$$
The second equality is a limit case of \eqtag{1.206}.
Wall polynomials are $q$-analogues of Laguerre polynomials on $[0,\infty)$,
so they might be called little $q$-Laguerre polynomials.

\sLP
3c)  Moak's \ref{Moa81}
{\sl $q$-Laguerre polynomials} (notation as in
\reff{GaRa90}{Exercise 7.43}) are given by
$$
\eqalignno{
L_n^\al(x;q)
:=&
{(q^{\al+1};q)_n\over(q;q)_n}\,
{}_1\phi_1(q^{-n};q^{\al+1};q,-xq^{n+\al+1})
\cr
=&
{1\over(q;q)_n}\,{}_2\phi_1(q^{-n},-x;0,q,q^{n+\al+1}).
\cr}
$$
The second equality is a limit case of \eqtag{1.204}.
They can be obtained from the Wall polynomials by replacing $q$ by $q^{-1}$.
Their orthogonality measure is not unique.
For instance, there is a continuous orthogonality measure
$$
\int_0^\infty L_m^\al(x;q)\,L_n^\al(x;q)\,
{x^\al\,dx\over(-(1-q)x;q)_\infty}=
\de_{m,n},\quad m\ne n,
$$
but also discrete orthogonality measures
$$
\int_0^\infty L_m^\al(cx;q)\,L_n^\al(cx;q)\,
{x^\al\,d_qx\over(-c(1-q)x;q)_\infty}=
\de_{m,n},\quad m\ne n,\quad c>0.
$$
Chihara \ref{Chi78}{Ch. VI, \S2} calls these polynomials
generalized Stieltjes-Wigert polynomials and he uses another notation.
For certain parameter values these polynomials may be considered
as $q$-analogues of the Charlier polynomials,
see \reff{GaRa90}{Exercise 7.13}.

\sLP
3d) {\sl Al-Salam-Carlitz I polynomials}
$$
U_n^{(a)}(x)=U_n^{(a)}(x;q):=
(-1)^n\,q^{n(n-1)/2}\,a^n\,{}_2\phi_1(q^{-n},x^{-1};0;q,qx/a)
$$
(cf.\ Al-Salam \& Carlitz \ref{AlCa65})
satisfy the three term recurrence relation
$$
x\,U_n^{(a)}(x)=U_{n+1}^{(a)}(x)+(1+a)\,q^n\,U_n^{(a)}(x)
-a\,q^{n-1}\,(1-q^n)\,U_{n-1}^{(a)}(x).
$$
Thus they are orthogonal polynomials for $a<0$.
They can be expressed in terms of big $q$-Jacobi polynomials by
$$
U_n^{(a)}(x)=\Pt_n(x;0,0,1,-a;q),
$$
so they are orthogonal \wrt the measure
$(qx,qx/a;q)_\infty d_qx$ on $[a,1]$, cf.\ \eqtag{2.198}.
By the above recurrence relation these polynomials are $q$-analogues
of Hermite polynomials of shifted argument,
so they may be considered as ``big'' $q$-Hermite polynomials.

\sLP
3e) {\sl Al-Salam-Carlitz II polynomials}
$$
V_n^{(a)}(x)=V_n^{(a)}(x;q):=U_n^{(a)}(x;q^{-1})
=(-1)^n\,q^{-n(n-1)/2}\,a^n\,{}_2\phi_0(q^{-n},x;-;q^na^{-1})
$$
(cf.\ \ref{AlCa65}).
For $a>0$ these form another family of $q$-analogues of the Charlier
polynomials.
On the other hand, the polynomials $x\mapsto V_n^{(\al/\bar \al)}(-q\al x)$
can be considered as $q$-analogues of Hermite polynomials of shifted
argument.
See Groenevelt \ref{Groe} for a further discussion of this case.

\sLP
4a) {\sl Stieltjes-Wigert polynomials}
$$
S_n(x;q):=(-1)^n\,q^{-n(2n+1)/2}\,
{}_1\phi_1(q^{-n};0;q,-q^{n+3/2}x)
$$
(cf.\ Chihara \reff{Chi78}{\S VI.2})
do not have a unique orthogonality measure.
It was already noted by Stieltjes that the corresponding moments
are an example of a non-determinate (Stieltjes) moment problem.
After suitable scaling these polynomials tend to Hermite polynomials as
$q\uparrow 1$.

\sLP
4b) {\sl Al-Salam-Carlitz I polynomials\/} $U_n^{(a)}$ with $a:=-1$
(these are
also known as {\sl discrete $q$-Hermite I polynomials} $h_n(\,.\,;q)$):
$$
\eqalignno{
h_n(x;q)=U_n^{(-1)}(x;q)
&=
q^{n(n-1)/2}\,{}_2\phi_1(q^{-n},x^{-1};0;q,-qx)
\cr
&=
\Pt_n(x;0,0,1,1;q)
\cr
&=
x^n\,{}_2\phi_0(q^{-n},q^{-n+1};-;q^2,q^{2n-1}x^{-2})
\cr}
$$
(cf.\ \ref{AlCa65}), $q$-analogues of Hermite polynomials.
By the last equality there is a quadratic transformation between these
polynomials and certain Wall polynomials. This is a limit case of
the quadratic transformations \eqtag{2.330} and \eqtag{2.340}.

\sLP
4c) {\sl Al-Salam-Carlitz II polynomials\/} $V_n^{(a)}$ of imaginary
argument and with $a:=-1$
(also known as {\sl discrete $q$-Hermite II polynomials}
$\tilde h_n(\,.\,;q)$):
$$
\eqalignno{
\tilde h_n(x;q)=i^{-n}\,V_n^{(-1)}(ix;q)
&=
i^{-n}\,q^{-n(n-1)/2}\,{}_2\phi_0(q^{-n},ix;-;-q^n)
\cr
&=\lim_{c\to\infty}\Pt_n(x;ic,ic,qc,qc;q)
\cr
&=x^n\,{}_2\phi_1(q^{-n},q^{-n+1};0;q^2,-x^{-2})
\cr}
$$
(cf.\ \ref{AlCa65}), also $q$-analogues of Hermite polynomials.
By the last equality there is a quadratic transformation between these
polynomials and certain of Moak's $q$-Laguerre polynomials.
This is a limit case of \eqtag{2.330} and \eqtag{2.340}.

\Subsec6 {The Askey-Wilson integral}
We remarked earlier that, whenever some nontrivial evaluation of an
integral can be given, the orthogonal polynomials having the integrand
as weight function may be worthwhile to study.
In this subsection we will give an evaluation of the integral which
corresponds to the Askey-Wilson polynomials.
After the original evaluation in
Askey \& Wilson \ref{AsWi85} several easier approaches were given,
cf.\ Gasper \& Rahman \reff{GaRa90}{Ch. 6} and the references
given there.
Our proof below borrowed ideas from Kalnins \& Miller \ref{KaMi89}
and Miller \ref{Mil89} but is still different from theirs.
Fix $0<q<1$.
Let
$$
w_{a,b,c,d}(z):=
{(z^2,z^{-2};q)_\infty
\over
(az,az^{-1},bz,bz^{-1},cz,cz^{-1},dz,dz^{-1};q)_\infty}\,.
\eqno\eq{2.480}
$$
We want to evaluate the integral
$$
I_{a,b,c,d}
:=
{1\over2\pi i}\,
\oint_{|z|=1} w_{a,b,c,d}(z)\,{dz\over z}\,,\quad
|a|,|b|,|c|,|d|<1.
\eqno\eq{2.490}
$$
Note that $I_{a,b,c,d}$ is analytic in the four complex variables $a,b,c,d$
when these are bounded in absolute value by 1.
It is also symmetric in $a,b,c,d$.

\Lemma{2.100}
$$
I_{a,b,c,d}=
{1-abcd\over(1-ab)\,(1-ac)\,(1-ad)}\,
I_{qa,b,c,d}.
\eqno\eq{2.500}
$$

\Proof
The integral
$$
\oint
{w_{q^{1/2}a,q^{1/2}b,q^{1/2}c,q^{1/2}d}(z)
\over
z-z^{-1}}\,
{dz\over z}
$$
equals on the one hand
$$
\displaylines{
\qquad\oint
{w_{q^{1/2}a,q^{1/2}b,q^{1/2}c,q^{1/2}d}(q^{1/2}z)
\over
q^{1/2}z-q^{-1/2}z^{-1}}\,
{dz\over z}\hfill
\cr
\hfill=\oint
{(1-az)\,(1-bz)\,(1-cz)\,(1-dz)\,w_{a,b,c,d}(z)
\over
q^{1/2}\,z\,(1-z^2)}\,
{dz\over z}\qquad
\cr}
$$
and on the other hand
$$
\displaylines{
\qquad\oint
{w_{q^{1/2}a,q^{1/2}b,q^{1/2}c,q^{1/2}d}(q^{-1/2}z)
\over
q^{-1/2}z-q^{1/2}z^{-1}}\,
{dz\over z}\hfill
\cr
\hfill=-\oint
{(1-a/z)(1-b/z)(1-c/z)(1-d/z)w_{a,b,c,d}(z)
\over
q^{1/2}\,z^{-1}\,(1-z^{-2})}\,
{dz\over z}\,.\qquad
\cr}
$$
Subtraction yields
$$
\displaylines{
0=
q^{-1/2}a^{-1}\,
\oint\bigl\{
-(1-abcd)(1-az)(1-az^{-1})
+(1-ab)(1-ac)(1-ad)\bigr\}\hfill
\cr
\hfill\times w_{a,b,c,d}(z)\,{dz\over z}
=\hfill2\pi iq^{-1/2}a^{-1}
\bigl\{-(1-abcd)I_{qa,b,c,d}
+(1-ab)(1-ac)(1-ad)I_{a,b,c,d}\bigr\}.
\hfill\halmos
\cr}
$$

\bPP
By iteration of \eqtag{2.500} and use of the analyticity and symmetry of
$I_{a,b,c,d}$ we obtain
$$
I_{a,b,c,d}
=
{(abcd;q)_\infty
\over
(ab,ac,ad,bc,bd,cd;q)_\infty}\,
I_{0,0,0,0}\,.
\eqno\eq{2.505}
$$
We might evaluate $I_{0,0,0,0}$ by the Jacobi triple product identity
\eqtag{1.300},
but it is easier to observe that
$I_{1,q^{1/2},-1,-q^{1/2}}=1$.
(Note that this case can be continuously reached from the domain of
definition of $I_{a,b,c,d}$ in \eqtag{2.490}.)
Hence \eqtag{2.505} yields that
$I_{0,0,0,0}=2/(q;q)_\infty$.
Thus
$$
I_{a,b,c,d}
=
{2(abcd;q)_\infty
\over
(ab,ac,ad,bc,bd,cd,q;q)_\infty}\,.
\eqno\eq{2.510}
$$

By the symmetry of $w_{a,b,c,d}(z)$ under $z\mapsto z^{-1}$
we can rewrite \eqtag{2.490}, \eqtag{2.510} as
$$
{1\over2\pi}\,
\int_0^\pi w_{a,b,c,d}(e^{i\th})\,d\th
=
{(abcd;q)_\infty
\over
(ab,ac,ad,bc,bd,cd,q;q)_\infty}\,,\quad
|a|,|b|,|c|,|d|<1.
\eqno\eq{2.520}
$$
Here $w$ is still defined by \eqtag{2.480}.
The integral \eqtag{2.520} is known as the {\sl Askey-Wilson integral}.

\Subsec7 {Askey-Wilson polynomials}
We now look for an orthogonal system corresponding to the
weight function in the Askey-Wilson integral.
First observe that
$$
\eqalignno{
{1\over 2\pi i}\,
\oint
{(az,az^{-1};q)_k\over(ac,ad;q)_k}\,
{(bz,bz^{-1};q)_l\over(bc,bd;q)_l}\,
w_{a,b,c,d}(z)\,
{dz\over z}
&=
{I_{q^ka,q^lb,c,d}
\over
(ac,ad;q)_k\,(bc,bd;q)_l}
\cr
&=
{(ab;q)_{k+l}\over(abcd;q)_{k+l}}\,
I_{a,b,c,d}.
\cr}
$$
Thus
$$
\eqalignno{
&{1\over I_{a,b,c,d}}\,
{1\over 2\pi i}\,
\oint
\biggl\{
\sum_{k=0}^n
{(q^{-n},q^{n-1}abcd;q)_k\,q^k
\over
(ab,q;q)_k}\,
{(az,az^{-1};q)_k\over(ac,ad;q)_k}
\biggr\}\,
{(bz,bz^{-1};q)_l\over(bc,bd;q)_l}\,
w_{a,b,c,d}(z)\,
{dz\over z}
\cr
&\qquad\qquad
=\sum_{k=0}^n
{(q^{-n},q^{n-1}abcd;q)_k\,q^k
\over
(ab,q;q)_k}\,
{(ab;q)_{k+l}\over(abcd;q)_{k+l}}
\cr
&\qquad\qquad
={(ab;q)_l\over(abcd;q)_l}\,
{}_3\phi_2\left[
{q^{-n},q^{n-1}abcd,q^lab\atop ab,q^labcd};q,q\right]
\cr
&\qquad\qquad
={(ab;q)_l\over(abcd;q)_l}\,
{(q^{-n+1}c^{-1}d^{-1},q^{-l};q)_n
\over
(ab,q^{-n-l+1}/(abcd);q)_n}
=0,\quad l=0,1,\ldots,n-1.
&\eq{2.525}
\cr}
$$
Here we used the $q$-Saalsch\"utz formula
\eqtag{1.400}.

The above orthogonality suggests to define
{\sl Askey-Wilson polynomials}
(Askey \& Wilson \ref{AsWi85})
$$
\eqalignno{
{p_n(\cos\th;a,b,c,d\mid q)
\over a^{-n}\,(ab,ac,ad;q)_n}
=&
r_n(\cos\th;a,b,c,d\mid q)
&\eq{2.529}
\cr
:=&
{}_4\phi_3\left[{
q^{-n},q^{n-1}abcd,ae^{i\th},ae^{-i\th}\atop ab,ac,ad};q,q\right]
&\eq{2.530}
\cr
=&\sum_{k=0}^n
{(q^{-n},q^{n-1}abcd;q)_k\,q^k
\over
(ab,q;q)_k}\,
{(ae^{i\th},ae^{-i\th};q)_k\over(ac,ad;q)_k}.
\cr}
$$
Since
$$
(ae^{i\th},ae^{-i\th};q)_k
=
\prod_{j=0}^{k-1}
(1-2aq^j\cos\th+a^2q^{2j}),
$$
formula \eqtag{2.530} defines a polynomial of degree $n$ in $\cos\th$.
It follows from \eqtag{2.525} that the functions
$\th\mapsto p_n(\cos\th;a,b,c,d\mid q)$ are orthogonal \wrt the measure
$w_{a,b,c,d}(e^{i\th})\,d\th$ on $[0,\pi]$, so we are really dealing with
orthogonal polynomials.
Now
$$
p_n(x;a,b,c,d\mid q)
=
k_n\,x^n+
\hbox{terms of lower degree,\quad with }k_n:=2^n\,(q^{n-1}abcd;q)_n,
\eqno\eq{2.540}
$$
so the coefficient of $x^n$ is symmetric in $a,b,c,d$.
Since the weight function
is also symmetric in $a,b,c,d$,
the Askey-Wilson polynomial will be itself symmetric in $a,b,c,d$.

\Prop{2.200}
Let $|a|,|b|,|c|,|d|<1$. Then
$${1\over2\pi}\,\int_0^\pi p_n(\cos\th)\,p_m(\cos\th)\,w(\cos\th)\,d\th=
\de_{m,n}\,h_n\,,$$
where
$$
\eqalignno{
p_n(\cos\th)&=p_n(\cos\th;a,b,c,d\mid q),\cr
\noalign{\vskip3pt}
w(\cos\th)&=
{(e^{2i\th},e^{-2i\th};q)_\infty
\over
(ae^{i\th},ae^{-i\th},be^{i\th},be^{-i\th},
ce^{i\th},ce^{-i\th},de^{i\th},de^{-i\th};q)_\infty}\,,\cr
{h_n\over h_0}&=
{(1-q^{n-1}abcd)\,(q,ab,ac,ad,bc,bd,cd;q)_n
\over(1-q^{2n-1}abcd)\,(abcd;q)_n}\,,\cr
\noalign{\hbox{and}}
h_0&={(abcd;q)_\infty\over
(q,ab,ac,ad,bc,bd,cd;q)_\infty}\,.
\cr}
$$
The orthogonality measure is positive if
$a,b,c,d$ are real, or if complex, appear in conjugate pairs,

\Proof
Apply \eqtag{2.525} to
$$
{1\over I_{a,b,c,d}}\,
{1\over 2\pi i}\,
\oint (p_np_m)\left((z+z^{-1})/2;a,b,c,d\mid q\right)\,
w_{a,b,c,d}(z)\,{dz\over z}\,,
$$
where $p_n$ is expanded according to \eqtag{2.530} and $p_m$ similarly,
but with $a$ and $b$ interchanged.\hhalmos

\bPP
Note that we can evaluate the Askey-Wilson polynomial
$p_n((z+z^{-1})/2)$ for $z=a$
(and by symmetry also for $z=b,c,d$):
$$
p_n((a+a^{-1})/2;a,b,c,d\mid q)=
a^{-n}\,(ab,ac,ad;q)_n.
\eqno\eq{2.550}
$$
When we write the three term recurrence relation for the Askey-Wilson polynomial
as
$$
2x\,p_n(x)=A_n\,p_{n+1}(x)+B_n\,p_n(x)+C_n\,p_{n-1}(x),
\eqno\eq{2.560}
$$
then the coefficients $A_n$ and $C_n$ can be computed from
$$
A_n={2k_n\over k_{n+1}},\quad
C_n={2k_{n-1}\over k_n}\,{h_n\over h_{n-1}},
$$
where $k_n$ and $h_n$ are given by \eqtag{2.540} and Proposition \thtag{2.200},
respectively.
Then $B_n$ can next be computed from \eqtag{2.550}
by substituting $x:=(a+a^{-1})/2$ in \eqtag{2.560}.

\Subsec8 {Various results}
Here we collect without proof some further results about Askey-Wilson
polynomials and their special cases and limit cases.

The {\sl $q$-ultraspherical polynomials}
$$
\eqalignno{
C_n(\cos\th;\be\mid q)
:=&
{(\be^2;q)_n\over\be^{n/2}(q;q)_n}\,
{}_4\phi_3\left[
{q^{-n},q^n\be^2,\be^{1/2}e^{i\th},\be^{1/2}e^{-i\th}
\atop
\be q^{1/2},-\be q^{1/2},-\be};q,q\right]
\cr
=&\const p_n(\cos\th;\be^{1/2},\be^{1/2}q^{1/2},-\be^{1/2},
-\be^{1/2}q^{1/2}\mid q)
\cr
=&
{(\be;q)_n\over(q;q)_n}\,
\sum_{k=0}^n{(q^{-n},\be;q)_k\over(q^{1-n}\be^{-1},q;q)_k}\,
(q/\be)^k\,e^{i(n-2k)\th}
\cr}
$$
are special Askey-Wilson polynomials which were already known to Rogers
(1894), however not as orthogonal polynomials.
See Askey \& Ismail \ref{AsIs83}.

By easy arguments using analytic continuation, contour deformation and taking of
residues it can be seen that Askey-Wilson polynomials for more general values
of the parameters than in Proposition \thtag{2.200} become orthogonal
\wrt a measure which contains both a continuous and discrete part
(Askey \& Wilson \ref{AsWi85})

\Prop{2.210}
Assume $a,b,c,d$ are real, or if complex, appear in conjugate pairs,
and that the pairwise products of $a,b,c,d$ are not $\ge1$, then
$$
{1\over2\pi}\,\int_0^\pi p_n(\cos\th)\,p_m(\cos\th)\,w(\cos\th)\,d\th+
\sum_k p_n(x_k)\,p_m(x_k)\,w_k=\de_{m,n}\,h_n\,,
$$
where
$p_n(\cos\th)$, $w(\cos\th)$ and $h_n$ are as in Proposition \thtag{2.200},
while the $x_k$ are the points $(eq^k+e^{-1}q^{-k})/2$ with $e$ any of the
parameters $a,b,c$ or $d$ whose absolute value is larger than one, the sum is
over the $k\in\Zplus$ with $|eq^k|>1$ and $w_k$ is $w_k(a;b,c,d)$ as defined by
\reff{AsWi85}{(2.10)} when $x_k=(aq^k+a^{-1}q^{-k})/2$.
(Be aware that $(1-aq^{2k})/(1-a)$ should be replaced by
$(1-a^2q^{2k})/(1-a^2)$ in \reff{AsWi85}{(2.10)}.)

\sPP
Both the big and the little $q$-Jacobi polynomials can be obtained as limit
cases of Askey-Wilson polynomials. In order to formulate these limits,
let $r_n$ be as in \eqtag{2.529}. Then
$$
\eqalignno{
\lim_{\la\to 0}
r_n&\left({q^{1/2}x\over2\la(cd)^{1/2}};
\la a(qd/c)^{1/2},\la^{-1}(qc/d)^{1/2},
-\la^{-1}(qd/c)^{1/2},-\la b(qc/d)^{1/2}\mid q\right)\cr
&\qquad=P_n(x;a,b,c,d;q)
\cr}
$$
and
$$
\lim_{\la\to 0}r_n\left({q^{1/2}x\over2\la^2};q^{1/2}\la^2a,q^{1/2}\la^{-2},
-q^{1/2},-q^{1/2}b\mid q\right)=
{(qb;q)_n\over(q^{-n}a^{-1};q)_n}\,p_n(x;b,a;q).
$$
See Koornwinder \reff{Koo93}{\S6}.
As $\la$ becomes smaller in these two limits, the number of mass points in the
orthogonality measure grows, while the support of the continuous part shrinks.
In the limit we have only infinitely many mass points and no continuous mass
left.

An important special class of Askey-Wilson polynomials are the
{\sl Al-Salam-Chihara polynomials\/} $p_n(x;a,b,0,0\mid q)$,
cf.\ \reff{AsIs84}{Ch. 3}.
Both these polynomials and the continuous $q$-ultraspherical polynomials
have the {\sl continuous $q$-Hermite polynomials}
$p_n(x;0,0,0,0\mid q)$
(cf.\ Askey \& Ismail \ref{AsIs83})
as a limit case.

The Askey-Wilson polynomials are eigenfunctions of a kind of $q$-difference
operator.
Write $P_n(e^{i\th})$ for the expression in \eqtag{2.530} and put
$$
A(z):=
{(1-az)(1-bz)(1-cz)(1-dz)\over(1-z^2)(1-qz^2)}\,.
$$
Then
$$
\displaylines{
\qquad A(z)P_n(qz)-(A(z)+A(z^{-1}))P_n(z)+A(z^{-1})P_n(q^{-1}z)\hfill
\cr
\hfill=-(1-q^{-n})(1-q^{n-1}abcd)P_n(z).\qquad
\cr}
$$
The operator
on the \LHS\ can be factorized as a product of two shift operators.
See Askey \& Wilson \reff{AsWi85}{\S5}.

An analytic continuation of the Askey-Wilson polynomials gives
{\sl $q$-Racah polynomials}
$$
R_n(q^{-x}+q^{x+1}\ga\de):=
{}_4\phi_3\left[{
q^{-n}q^{n+1}\al\be,q^{-x},q^{x+1}\ga\de\atop
\al q,\be\de q,\ga q};q,q\right],
$$
where one of $\al q$, $\be\de q$ or $\ga q$ is $q^{-N}$ for some $N\in\Zplus$
and $n=0,1,\ldots,N$ (see Askey \& Wilson \ref{AsWi79}).
These satisfy an orthogonality of the form
$$
\sum_{x=0}^N R_n(\mu(x))\,R_m(\mu(x))\,w(x)=\de_{m,n}\,h_n,\quad
m,n=0,1,\ldots,N,
$$
where $\mu(x):=q^{-x}+q^{x+1}\ga\de$.
They are eigenfunctions of a second order difference operator in $x$.

There is nice characterization theorem of Leonard
\ref{Leo82}
for the $q$-Racah polynomials.
Let $N\in\Zplus$ or $N=\infty$. Let the polynomials $p_n$ ($n\in\Zplus$,
$n<N+1$) be orthogonal \wrt weights on distinct points $\mu_k$
($k\in\Zplus$, $k<N+1$).
Let the polynomials $p_n^*$ be similarly orthogonal \wrt weights on distinct
points $\mu_k^*$.
Suppose that the two systems are dual in the sense that
$$
p_n(\mu_k)=p_k^*(\mu_n^*).
$$
Then the $p_n$ are $q$-Racah polynomials or one of their limit cases.

\bPP\goodbreak\noindent \exnumber1
{\bf Exercises to \S\the\sectionnumber}
\par\nobreak\smallskip\noindent\the\sectionnumber.\the\exnumber\quad
\advance\exnumber by 1
Show that (at least formally) the generating function
$$
e^{2xt-t^2}=\sum_{n=0}^\infty{H_n(x)\over n!}\,t^n
$$
for Hermite polynomials (cf.\ \reff{Erd2}{10.13(19)})
follows from the generating function
$$
(1-t)^{-\al-1}\,e^{-tx/(1-t)}=\sum_{n=0}^\infty L_n^\al(x)\,t^n,\quad |t|<1,
$$
for Laguerre polynomials (cf.\ \reff{Erd2}{10.12(17)})
by the limit transition
$$
H_n(x)=(-1)^n\,2^{n/2}\,n!\,\lim_{\al\to\infty}
\al^{-n/2}\,L_n^\al((2\al)^{1/2}x+\al)
$$
given in \eqtag{2.60}.

\nextex
Show that (at least formally) the above generating function
for Laguerre polynomials follows from the generating function
for Jacobi polynomials
$$
\displaylines{
\hfill\sum_{n=0}^\infty P_n^{(\al,\be)}(x)\,t^n=
2^{\al+\be}\,R^{-1}\,(1-t+R)^{-\al}\,(1+t+R)^{-\be},\quad |t|<1,\hfill
\cr
\hbox{where}\hfill
\cr
\hfill R:=(1-2xt+t^2)^{1/2},\hfill
\cr}
$$
(cf.\ \reff{Erd2}{10.8 (29)}) by the limit transition
$$
L_n^\al(x)=\lim_{\be\to\infty}P_n^{(\al,\be)}(1-2x/\be)
$$
given in \eqtag{2.40}.

\nextex
Show that (at least formally) the above generating function
for Hermite polynomials follows from the specialization $\al=\be$ of the above
generating function for Jacobi polynomials by the limit transition
$$
H_n(x)=
2^n\,n!\,\lim_{\al\to\infty}\al^{-n/2}\,P_n^{(\al,\al)}(\al^{-1/2}x)
$$
given in \eqtag{2.50}.

\nextex
Prove that the Charlier polynomials
$$
C_n(x;a):={}_2F_0(-n,-x;-;-a^{-1}),\quad a>0,
$$
satisfy the recurrence relation
$$
x\,C_n(x;a)=-a\,C_{n+1}(x;a)+(n+a)\,C_n(x;a)-n\,C_{n-1}(x;a).
$$

\nextex
Prove that
$$
\lim_{a\to\infty}(-(2a)^{1/2})^n\,C_n((2a)^{1/2}x+a)=H_n(x)
$$
by using the above recurrence relation for Charlier polynomials and the
recurrence relation
$$
H_{n+1}(x)-2x\,H_n(x)+2n\,H_{n-1}(x)=0
$$
for Hermite polynomials.

\nextex
Prove the following generating function for Al-Salam-Chihara polynomials
(notation as for Askey-Wilson polynomials in \eqtag{2.529}, \eqtag{2.530}):
$$
{(zc,zd;q)_\infty\over(ze^{i\th},ze^{-i\th};q)_\infty}=
\sum_{m=0}^\infty z^m\,
{p_m(\cos\th;0,0,c,d\mid q)\over(q;q)_m}\,,\quad |z|<1.
$$

\nextex
Use the above generating function in order to
derive the following transformation formula from little $q$-Jacobi
polynomials (notation and definition by \eqtag{2.300}) to Askey-Wilson
polynomials by means of a summation kernel involving
Al-Salam-Chihara polynomials:
$$
\displaylines{
\qquad
{a^n\,p_n(\cos\th;a,b,c,d\mid q)\over(ab,ac,ad;q)_n}\,
{(ac,ad;q)_\infty\over(ae^{i\th},ae^{-i\th};q)_\infty}\hfill
\cr
\hfill
=\sum_{m=0}^\infty p_n(q^m;q^{-1}ab,q^{-1}cd;q)\,
a^m\,{p_m(\cos\th;0,0,c,d\mid q)\over(q;q)_m}\,.\qquad
\cr}
$$

{\Ref
\refitem{AlCa65}
W. A. Al-Salam \&
L. Carlitz,
{\sl Some orthogonal $q$-polynomials},
Math. Nachr. 30 (1965), 47--61.

\refitem{And86}
G. E. Andrews,
{\sl $q$-Series: their development and application in analysis, number theory,
combinatorics, physics, and computer algebra},
Regional Conference Series in Math. 66,
Amer. Math. Soc., 1986.

\refitem{AnAs77}
G. E. Andrews \&
R. Askey,
{\sl Enumeration of partitions: The role of Eulerian series and $q$-orthogonal
polynomials},
in {\sl Higher combinatorics},
M. Aigner (ed.),
Reidel, 1977, pp. 3--26.

\refitem{AnAs85}
G. E. Andrews \&
R. Askey,
{\sl Classical orthogonal polynomials},
in {\sl Polyn{\^o}mes orthogonaux et applications},
C. Brezinski,
A. Draux,
A. P. Magnus,
P. Maroni \&
A. Ronveaux (eds.),
Lecture Notes in Math.
1171,
Springer-Verlag,
1985,
pp. 36--62.

\refitem{AsIs83}
R. Askey \&
M. E. H. Ismail,
{\sl A generalization of ultraspherical polynomials},
in {\sl Studies in Pure Mathematics},
P. Erd{\"o}s (ed.),
Birkh{\"a}user, 1983, 55--78.

\refitem{AsIs84}
R. Askey \&
M. E. H. Ismail,
{\sl Recurrence relations, continued fractions and orthogonal polynomials},
Mem. Amer. Math. Soc. 49 (1984), no. 300.

\refitem{AsWi79}
R. Askey \&
J. Wilson,
{\sl A set of orthogonal polynomials that generalize the Racah coefficients or
6-j symbols},
SIAM J. Math. Anal. 10 (1979), 1008--1016.

\refitem{AsWi85}
R. Askey \&
J. Wilson,
{\sl Some basic hypergeometric orthogonal polynomials that generalize Jacobi
polynomials},
Mem. Amer. Math. Soc. 54 (1985), no. 319.

\refitem{AtSu85}
N. M. Atakishiyev \&
S. K. Suslov,
{\sl The Hahn and Meixner polynomials of an imaginary argument and some of their
applications},
J. Phys. A 18 (1985), 1583--1596.

\refitem{AtRaSu93}
N. M. Atakishiyev, M. Rahman \&
S. K. Suslov,
{\sl On classical orthogonal polynomials},
Constr. Approx. 11 (1995), 181--226.

\refitem{Bai35}
W. N. Bailey,
{\sl Generalized hypergeometric series},
Cambridge University Press, 1935;
reprinted by Hafner Publishing Company, 1972.

\refitem{Chi78}
T. S. Chihara,
{\sl An introduction to orthogonal polynomials},
Gordon and Breach, 1978;
reprinted by Dover Publications, 2011.

\refitem{Erd1}
A. Erd{\'e}lyi,
W. Magnus,
F. Oberhettinger \&
F. G. Tricomi,
Higher transcendental functions, Vol. 1,
McGraw-Hill, 1953.

\refitem{Erd2}
A. Erd{\'e}lyi,
W. Magnus,
F. Oberhettinger \&
F. G. Tricomi,
Higher transcendental functions, Vol. 2,
McGraw-Hill, 1953.

\refitem{GaRa90}
G. Gasper \&
M. Rahman,
{\sl Basic hypergeometric series},
Encyclopedia of Mathematics and its Applications 35,
Cambridge University Press, 1990; second edition 2004.

\refitem{Hah49}
W. Hahn,
{\sl {\"U}ber Orthogonalpolynome, die $q$-Differenzengleichungen gen{\"u}gen},
Math. Nachr. 2 (1949), 4--34, 379.

\refitem{Ism77}
M. E. H. Ismail,
{\sl A simple proof of Ramanujan's ${}_1\psi_1$ sum},
Proc. Amer. Math. Soc. 63 (1977), 185--186.

\refitem{KaMi89}
E. G. Kalnins \&
W. Miller, Jr.,
{\sl Symmetry techniques for $q$-series: Askey-Wilson polynomials},
Rocky Mountain J. Math. 19 (1989), 223--230.

\refitem{KoeSwa}
R. Koekoek \& R. F. Swarttouw,
{\sl The Askey scheme of hypergeometric orthogonal polynomials and its
$q$-analogue},
Report 94-05, Delft University of Technology,
Faculty of Technical Mathematics
and Informatics, 1994;
an extended version appeared as Report 98-17, 1998.
A further extension is in the book
R. Koekoek, P. A. Lesky \& R. F. Swarttouw,
{\sl Hypergeometric orthogonal polynomials and their $q$-analogues},
Springer-Verlag, 2010.

\refitem{Koo88}
T. H. Koornwinder,
{\sl Group theoretic interpretations of Askey's scheme of hypergeometric
orthogonal polynomials},
in {\sl Orthogonal polynomials and their applications},
M. Alfaro,
J. S. Dehesa,
F. J. Marcellan,
J. L. Rubio de Francia \&
J. Vinuesa (eds.),
Lecture Notes in Math. 1329, Springer-Verlag, 1988, pp. 46--72.

\refitem{Koo89}
T. H. Koornwinder,
{\sl Representations of the twisted $SU(2)$ quantum group and some
$q$-hyper\-geo\-metric orthogonal polynomials},
Indag. Math. 51 (1989), 97--117.

\refitem{Koo90}
T. H. Koornwinder,
{\sl Jacobi functions as limit cases of $q$-ultraspherical polynomials},
J.~Math. Anal. Appl. 148 (1990), 44--54.

\refitem{Koo93}
T. H. Koornwinder,
{\sl Askey-Wilson polynomials as zonal spherical functions on the $SU(2)$
quantum group},
SIAM J. Math. Anal.
24
(1993),
795--813.

\refitem{KoSw}
T. H. Koornwinder \&
R. F. Swarttouw,
{\sl On $q$-Analogues of the Fourier and Hankel transforms},
Trans. Amer. Math. Soc.
333
(1992),
445--461; corrected version appeared as
arXiv:1208.2521v1 [math.CA], 2012.

\refitem{Lab90}
J. Labelle,
{\sl Askey's scheme of hypergeometric orthogonal polynomials},
poster, Universit\'e de Quebec \`a Montr\'eal, 1990.

\refitem{Leo82}
D. A. Leonard,
{\sl Orthogonal polynomials, duality and association schemes},
SIAM J. Math. Anal. 13 (1982), 656--663.

\refitem{Mil89}
W. Miller, Jr.,
{\sl Symmetry techniques and orthogonality for $q$-Series},
in {\sl $q$-Series and partitions},
D. Stanton (ed.),
IMA Volumes in Math. and its Appl. 18,
Springer-Verlag, 1989, pp. 191--212.

\refitem{Mim89}
K. Mimachi,
{\sl Connection problem in holonomic $q$-difference system associated with
a Jackson integral of Jordan-Pochhammer type},
Nagoya Math. J. 116 (1989), 149--161.

\refitem{Moa81}
D. S. Moak,
{\sl The $q$-analogue of the Laguerre polynomials},
J. Math. Anal. Appl. 81 (1981), 20--47.

\refitem{Olv74}
F. W. J. Olver,
{\sl Asymptotics and special functions},
Academic Press, 1974.

\refitem{RaVe}
M. Rahman \&
A. Verma,
{\sl Quadratic transformation formulas for basic hypergeometric series},
Trans. Amer. Math. Soc. 335 (1993), 277--302.

\refitem{Tit39}
E. C. Titchmarsh,
{\sl The theory of functions},
Oxford University Press, second ed., 1939.

\refitem{Wil80}
J. A. Wilson,
{\sl Some hypergeometric orthogonal polynomials},
SIAM J. Math. Anal. 11 (1980), 690--701.

\refitem{Groe}
W. Groenevelt,
{\sl Orthogonality relations for Al-Salam-Carlitz polynomials of type II},
arXiv:1309.7569v1 [math.CA], 2013.

\par
}
\bye